\documentclass[12pt,reqno]{amsart}
\usepackage{amsmath, amssymb, amsthm}
\usepackage{mathtools}
\usepackage{enumitem}
\usepackage[breaklinks]{hyperref}
\usepackage{xcolor}
\newcommand{\Gal}{\operatorname{Gal}}

\newcommand{\g}{\gamma}

\let\amstextb\b
\renewcommand{\b}{\ifmmode\beta\else\expandafter\amstextb\fi}
\newcommand{\ZZ}{\mathbb{Z}}
\newcommand{\Sel}{\operatorname{Sel}}
\newcommand{\QQ}{\mathbb{Q}}
\newcommand{\FF}{\mathbb{F}}
\newcommand{\SL}{\operatorname{SL}}
\newcommand{\x}{\xi}
\newcommand{\z}{\zeta}
\newcommand{\HCal}{{\mathcal{H}}}
\newcommand{\Ob}{\mathcal{O}}
\newcommand{\Diff}{\mathfrak{d}}
\newcommand{\nf}{\mathfrak{n}}
\newcommand{\fa}{\mathfrak{a}}
\newcommand{\fb}{\mathfrak{b}}
\newcommand{\pf}{\mathfrak{p}}
\newcommand{\qf}{\mathfrak{q}}
\newcommand{\mf}{\mathfrak{m}}
\newcommand{\lf}{\mathfrak{l}}
\newcommand{\Qf}{\mathfrak{Q}}
\newcommand{\Tr}{\operatorname{Tr}}
\newcommand{\Norm}{\mathrm{N}}
\newcommand{\ub}[1]{u(#1)}
\newcommand{\Gz}[1]{\Gamma_0(#1)}
\newcommand{\Mk}[2]{M_k\!\left(\Gamma_0(#1),\,#2\right)}
\newcommand{\ord}{\operatorname{ord}}
\newcommand{\ba}{{\mathbf a}}

\newcommand{\vp}{\varpi}                 
\newcommand{\kp}{\kappa_{\pi}}           
\newcommand{\Dv}{D_{\vp}}               
\newcommand{\ko}{\kappa_{0}}            
\newcommand{\Mo}{M_{0}}
\newcommand{\Kx}[1]{K_{#1}}              
\newcommand{\etwo}{e_{2}}
\newcommand{\idx}{\iota}
\newcommand{\ew}{\mathbf{e}_{w}}        
\newcommand{\tq}{q}                      
\newcommand{\aR}{a^{R}}
\newcommand{\Gs}{\mathbf{G}}
\newcommand{\Es}{\mathbf{E}}              

\newcommand{\nAdd}{\cite[Lemma~2.3]{n}}
\newcommand{\nPhase}{\cite[Lemma~2.6]{n}}

\numberwithin{equation}{section}
\theoremstyle{plain}
\newtheorem*{maintheorem}{Main Theorem}
\newtheorem*{bhcriterion}{Theorem}
\newtheorem{theorem}{Theorem}[section]
\newtheorem{lemma}[theorem]{Lemma}
\newtheorem{corollary}[theorem]{Corollary}
\newtheorem{proposition}[theorem]{Proposition}
\theoremstyle{definition}
\newtheorem{definition}[theorem]{Definition}
\newtheorem{example}[theorem]{Example}
\theoremstyle{remark}

\theoremstyle{plain}
\newtheorem*{conjecture}{Conjecture}
\title[Relative class numbers of CM extensions]
{Relative class numbers of CM extensions:
	valuation bounds and indivisibility}
\date{}
\hypersetup{colorlinks=false,pdfborder={0 0 1},
	citebordercolor={0 1 0},linkbordercolor={1 0 0}}
\author{T\`{u}ng Ho\`{a}ng}

\begin{document}
		\begin{abstract}
		
			Let $F$ be a totally real field of degree $g\ge2$. For a wide range of prescribed local conditions, a nonzero witness sum $w$ of Eisenstein coefficients gives at least $C\sqrt X/\log X$ CM quadratic extensions $K/F$ satisfying these conditions, with relative discriminant norm below $X$ and $v_p(h^-(K/F))\le v_p(w)$. This holds for every odd prime $p$ and every $X\ge X_0$, with $C$ and $X_0$ independent of $p$.
			
			We strengthen Takai's result by removing the assumption that $F/\QQ$ is Galois and improving the counting exponent from $\frac{1}{2g}$ to $\frac{1}{2}$, matching that of the classical Kohnen--Ono theorem. We modify Takai's argument by combining the operators $U_\pi$ and $V_\pi$ with a restriction direction adapted to the auxiliary prime element $\pi$.
			
			Kohnen and Ono proved the indivisibility bound for imaginary quadratic fields for each fixed prime $p>3$. Using witnesses computed in PARI/GP, we extend their bound to CM quadratic extensions of more than $20{,}000$ real quadratic fields. Over each of these fields, the bound holds for every odd prime $p$.
			
			Under suitable hypotheses, the same bound holds for quadratic twists with trivial $\ell$-Selmer group, with explicit examples up to degree $12$. For every odd prime $p$ and $n\ge1$, we prove the class-number indivisibility bound at $p$ over infinitely many cyclic fields of degree $p^n$ without a witness hypothesis. For $p\ge5$, the Selmer bound also holds over these fields under local and torsion hypotheses. For the fixed curve $123.a1$ over $\mathbb Q$, the $5$-Selmer bound holds over infinitely many cyclic fields of each degree $5^n$.

	\end{abstract}
	\maketitle
	
	\section{Introduction}

	The study of class numbers of imaginary quadratic fields has its origins in Gauss's theory of binary quadratic forms. Gauss conjectured that $h(D)\to\infty$ as $D\to-\infty$, where $D$ ranges over negative fundamental discriminants and $h(D)$ denotes the class number of $\mathbb Q(\sqrt D)$. Equivalently, only finitely many imaginary quadratic fields have any given class number. Heilbronn~\cite{Heilbronn} proved this conjecture in 1934.
	
	One can also ask how often $h(D)$ is coprime to a fixed odd prime $p$. When imaginary quadratic fields are ordered by absolute discriminant, the Cohen--Lenstra heuristics~\cite{CohenLenstra} predict that the proportion with this property is
	\[
	\prod_{n=1}^{\infty}(1-p^{-n})>0.
	\]
	This motivates the problem of finding lower bounds for the number of such fields with absolute discriminant less than $X$.
	
	An important result toward the Cohen--Lenstra prediction was obtained by Kohnen and Ono~\cite{KohnenOno}, who proved that for every prime $p>3$,
	$$
	\#\left\{D<0:\ |D|<X,\ p\nmid h(D)\right\}
	\gg_p \frac{\sqrt X}{\log X}.
	$$
	Their method uses modular forms of half-integral weight, the two operators $U_\ell$ and $V_\ell$, and Sturm's theorem. An explicit formula for the coefficient of $q^\ell$ in $(U_\ell-V_\ell)\theta^3$ allows them to choose $\ell$ so that this coefficient is a $p$-unit.
	
	A more general setting is obtained by replacing $\QQ$ with a totally real number field $F$. For a CM quadratic extension $K/F$, the corresponding invariant is the relative class number $h^-(K/F)=h(K)/h(F)$.
	
	Takai~\cite[Corollary~1]{Takai} generalized Kohnen--Ono's method to CM extensions $K/F$ when $F/\QQ$ is Galois. He used Hilbert modular forms of half-integral weight and diagonal restriction to apply Sturm's bound as in the classical case over $\QQ$. He proved that, for every sufficiently large prime $p$,
	\[
	\#\left\{
	K/F\text{ CM quadratic}:
	\mathrm{N}_{F/\QQ}(D(K/F))<X,\ 
	p\nmid h^-(K/F)
	\right\}
	\gg_{F,p}\frac{X^{\frac{1}{2g}}}{\log X},
	\]
	where $g=[F:\QQ]$.
	We improve Takai's exponent from $\frac{1}{2g}$ to $\frac{1}{2}$ and remove the assumption that $F/\QQ$ is Galois. For each odd prime $p$, we count extensions satisfying the prescribed local conditions and $v_p(h^-(K/F))\le v_p(w)$, where $w$ is a nonzero witness. This gives indivisibility when $v_p(w)=0$ and still bounds the valuation when $p$ divides $w$.
	
	\medskip
	\noindent\textbf{Main result.}
	
	Fix $F$, $\mathcal E\in\{\Es,\Gs\}$, and a local condition $(\nf,R)$ as in \S\ref{ss:sieve}. We use condition \textup{(W)} of \S\ref{sec:main}, with data $m_0,\Qf$ and nonzero witness sum $w$ given by \eqref{eq:W}. For an odd prime $p$ and an integer $t\ge0$, let $N_{p,t}(X)$ count the distinct fields $K=F(\sqrt{-\eta})$ with $\eta\in\Ob_+\cap R$, $\Norm D(K/F)<X$, and $v_p(h^-(K/F))\le t$.
	
	\begin{maintheorem}
		Assume \textup{(W)}, and fix $m_0,\Qf,w$ as above. Then $w\in\ZZ[1/2]$. There are constants $C>0$ and $X_0>1$, depending on $F,\mathcal E,(\nf,R),m_0,\Qf$, but not on $p$, such that
		\[
		N_{p,v_p(w)}(X)\ge C\frac{\sqrt X}{\log X}
		\qquad\text{for every odd prime $p$ and every $X\ge X_0$.}
		\]
		In particular, if $v_p(w)=0$, then
		\[
		\#\left\{K=F(\sqrt{-\eta})\ \middle|\
		\begin{aligned}
			&\eta\in\Ob_+\cap R,\quad \Norm D(K/F)<X,\\
			&p\nmid h^-(K/F)
		\end{aligned}\right\}
		\ge C\frac{\sqrt X}{\log X}\qquad(X\ge X_0).
		\]
	\end{maintheorem}
	
	Under the stated nonvanishing hypothesis, this extends the Kohnen--Ono method to every totally real field of degree at least two, with the same counting exponent as over $\QQ$. We choose the diagonal restriction to depend on the auxiliary prime element $\pi$. If $\Norm(\pi)=\ell$, the fourth power of the restricted form is modular on $\Gamma_0(\ell\Mo(\mf))$. Sturm's bound and a trace estimate then give a CM quadratic extension whose relative discriminant has norm $O_{F,\mf}(\ell^2)$. Each field arises from at most two auxiliary primes, so Chebotarev's theorem gives the lower bound $\sqrt X/\log X$.

	The nonvanishing hypothesis can be checked at a single coefficient. Lemma~\ref{lem:witness} shows that a nonzero witness sum exists if and only if $a_{\eta/2}\ne0$ for some $\eta\in\Ob_+\cap R$. This criterion does not involve reduction modulo $p$. Together with the Main Theorem, this gives the indivisibility bound for all but finitely many odd primes.
	
	Even when $p$ divides $w$, we still obtain the bound $v_p(h^-(K/F))\le v_p(w)$ for at least $C\sqrt X/\log X$ fields. The extension and Frobenius condition used in Chebotarev's theorem do not depend on $p$. For each $p$, we exclude at most $g$ auxiliary primes lying above $p$. Thus the constant and the threshold for $X$ can be chosen independently of $p$, although the fields counted may depend on $p$.
	
	Instead of requiring a single witness sum to be a $p$-unit for every odd prime $p$, we use several witness sums. If $w_1,\ldots,w_r$ are nonzero witnesses for the same local conditions, then the indivisibility bound holds whenever $v_p(w_i)=0$ for at least one $i$. Thus an odd prime $p$ divides the gcd of the numerators of the $w_i$ in lowest terms if and only if none of these witnesses is a $p$-unit. If this gcd is $2^k$, the bound holds for every odd prime $p$, even when no individual witness has this property. Since only finitely many witnesses are used, both the constant $C$ and the threshold for $X$ can be chosen independently of $p$.
	
	Fix the field, the Eisenstein series, and the local conditions, and assume \textup{(W)}. The admissible witnesses generate an ideal $I_w\subseteq\ZZ[1/2]$ as in Definition~\ref{def:witness-ideal}. Write $I_w=d_w\ZZ[1/2]$, where $d_w$ is positive and odd. Corollary~\ref{cor:witness-gcd} gives the counting bound with $v_p(h^-(K/F))\le v_p(d_w)$, the smallest valuation bound obtained from all admissible witnesses. If $p\nmid d_w$, this gives indivisibility. If $p\mid d_w$, the witnesses do not settle indivisibility at $p$. In particular, if finitely many witnesses generate $\ZZ[1/2]$, the indivisibility bound holds for every odd prime.	
	
	Using PARI/GP, we show that the Kohnen--Ono indivisibility  bound extends to CM quadratic extensions of more than  $20{,}000$ real quadratic fields. The code is available at \url{https://tunghoangmath.github.io/code.html}.
	
	\begin{theorem}\label{thm:real-quadratic-computation}
		For more than $20{,}000$ real quadratic fields $F$, including
		$\QQ(\sqrt2)$, there are constants $C_F>0$ and $X_F>1$ such that
		\[
		\#\left\{K/F\text{ CM quadratic}:
		\Norm_{F/\QQ}D(K/F)<X,\quad p\nmid h^-(K/F)\right\}
		\ge C_F\frac{\sqrt X}{\log X}
		\]
		for every odd prime $p$ and every $X\ge X_F$.
		
		More precisely, let $w_1,\ldots,w_r$ be witnesses for the fixed sieve
		whose numerators in lowest terms have gcd a power of $2$. Then one may take
		\[
		C_F=\frac{1}{8\ko\max_{1\le i\le r}[L_i:F]},
		\]
		where $\ko$ is the cutoff for this sieve and $L_i/F$ is the extension
		constructed in \S\ref{ss:proof-main} from the data defining $w_i$.
	\end{theorem}
	\begin{conjecture}
		For every real quadratic field $F$, there are constants $C_F>0$ and $X_F>1$ such that the indivisibility bound in Theorem~\ref{thm:real-quadratic-computation} holds for every odd prime $p$ and every $X\ge X_F$.
	\end{conjecture}
	Example~\ref{ex:sqrt2-witness-gcd} gives the details for $\QQ(\sqrt2)$.
	In this case, the two witnesses $6$ and $10$ suffice, since every odd
	prime is coprime to at least one of them.
	
	In \S\ref{sec:prime-power}, we combine Wiles's theorem with the Galois action on Fourier indices to construct $p$-unit witnesses. Let $p$ be an odd prime and let $F/\QQ$ be cyclic of degree $p^n$, with $n\ge1$. Suppose that exactly one rational prime $r$ ramifies in $F$, with $r=p$ or $r\equiv1\pmod4$. We obtain the lower bound $\sqrt X/\log X$ for the number of CM quadratic extensions $K/F$ with $\Norm D(K/F)<X$ and $p\nmid h(K)$, even with suitable prescribed local conditions (Theorem~\ref{thm:prime-power-local-counting}).
	
	For the quadratic twists of elliptic curves in \cite[Theorem~4.1]{n}, the Main Theorem improves the lower bound for the number of twists with trivial $\ell$-Selmer group from $\frac{X^{\frac{1}{2g}}}{\log X}$ to $\frac{\sqrt{X}}{\log X}$, where $g=[F:\QQ]$.
	
	Let $\pf_1,\ldots,\pf_r$ be the distinct prime ideals
	dividing $\ell N(E)$, and let $\chi_i$ and $\varepsilon_i$
	be the corresponding characters and signs chosen as in
	\cite[Theorem~4.1]{n}. Writing $\mathfrak f_i$ for the
	conductor of $\chi_i$, put
	\begin{equation}\label{eq:selmer-local-condition}
		\nf=\prod_{i=1}^r\mathfrak f_i,\qquad
		R=\{\eta\in(\Ob/\nf)^\times:
		\chi_i(\eta)=\varepsilon_i\ \text{for all }i\}.
	\end{equation}
	As verified in \S\ref{ss:selmer-counting}, the pair
	$(\nf,R)$ is a local condition in the sense of
	Definition~\ref{def:local}.
	
	\begin{theorem}\label{thm:selmer-twists}
		Let $F$ be a totally real Galois number field of degree
		$g\ge2$, and let $\ell\ge5$ be a prime such that
		$\ell\nmid h(F)$. Let $E/F$ be an elliptic curve with
		an $F$-rational point of order $\ell$, satisfying
		conditions \textup{(i)} and \textup{(ii)} of
		\cite[Theorem~3.6]{n}.
		Let $(\nf,R)$ be the local condition defined in
		\eqref{eq:selmer-local-condition}.
		Assume that \textup{(W)} holds for one of $\Es$ and
		$\Gs$, with $v_\ell(w)=0$. Then
		\[
		\#\left\{
		d\in F^\times/(F^\times)^2:
		\begin{array}{l}
			d\text{ totally negative},\\
			\Norm_{F/\QQ}D(F(\sqrt d)/F)<X,\\
			\Sel_\ell(E^d,F)=0
		\end{array}
		\right\}
		\gg_{F,E,\ell}\frac{\sqrt X}{\log X}.
		\]
	\end{theorem}

	There is no upper bound on $g$ in this theorem. In
	\S\ref{ss:selmer-witness-search}, we give a PARI/GP algorithm for
	searching for $\ell$-unit witnesses using inert rational primes.
	Table~\ref{tab:selmer-witnesses} gives examples for six elliptic curves
	over five totally real fields, including fields of degrees $6$ and
	$12$, with $\ell=13$ and $\ell=37$, respectively.
	
	For certain cyclic fields of degree $p^n$, we obtain the same
	Selmer bound at $\ell=p$ without assuming \textup{(W)}, under the
	local hypotheses of Theorem~\ref{thm:prime-power-selmer}.
	In particular, for the fixed elliptic curve $E/\QQ$ with LMFDB label
	\texttt{123.a1}, Corollary~\ref{cor:prime-power-five-selmer} gives the
	bound $\sqrt X/\log X$ for twists with trivial $5$-Selmer group over
	every nontrivial finite layer of the cyclotomic $\ZZ_5$-extension of
	$\QQ$ and over infinitely many distinct cyclic fields of each degree
	$5^n$. No separate witness computation is needed for these fields.
	
	The rest of this paper is organized as follows.
	
	Section~\ref{sec:prelim} introduces the Eisenstein series and the local conditions used throughout the paper. Section~\ref{sec:operators} studies the two operators $U_\pi$ and $V_\pi$ and their action on Fourier coefficients. Section~\ref{sec:restriction} combines twisted diagonal restriction with Sturm's bound to produce coefficients with controlled $p$-adic valuation. Section~\ref{sec:main} proves the nonvanishing criterion and the main counting theorem, and explains how to choose and combine witnesses, with an explicit example over $\QQ(\sqrt2)$. Section~\ref{sec:prime-power} constructs $p$-unit witnesses over cyclic fields of prime-power degree. Finally, Section~\ref{sec:selmer} gives the applications to Selmer groups.
	
	\section{Preliminaries}\label{sec:prelim}
	
	\subsection{Notation}\label{ss:notation}
	Throughout, $F$ is a totally real number field of degree $g=[F:\QQ]\geq2$, with ring of integers $\Ob$, different $\Diff$ and discriminant $D(F)$; $\ba$ is the set of real places of $F$, and $\Ob_{+}$ is the set of totally positive elements of $\Ob$. For $x\in F$ we write $e(x)=\exp(2\pi i\Tr(x))$, and $e(\x z)=\exp(2\pi i\sum_{v\in\ba}\x_vz_v)$ for $z\in\HCal^{\ba}$. For an integral ideal $\mf$ with $4\Ob\mid\mf$,
	\begin{equation*}
		\Gz{\mf}=\Big\{\begin{psmallmatrix}a&b\\c&d\end{psmallmatrix}\in\SL_2(F):
		a,d\in\Ob,\ b\in 2\Diff^{-1},\ c\in 2^{-1}\mf\Diff\Big\}.
	\end{equation*}
	The theta series is
	\begin{equation*}
		\theta(z)=\sum_{x\in\Ob}e(x^2z/2),
	\end{equation*}
	and we put $j(\g,z)=\prod_{v\in\ba}(c_vz_v+d_v)$. For $\g\in\Gz{4\Ob}$, let $h(\g,z)$ be the factor of automorphy supplied by \cite[Proposition~1.2]{Sh}. It satisfies
	\begin{equation*}
		\theta(\g z)=h(\g,z)\theta(z)
	\end{equation*}
	and
	\begin{equation}\label{eq:h2}
		h(\g,z)^2=\omega(d_\g)\,j(\g,z)\qquad(\g\in\Gz{4\Ob}),
	\end{equation}
	with a sign $\omega(d_\g)\in\{\pm1\}$ which depends only on $d_\g$. It also satisfies the cocycle relation
	\begin{equation*}
		h(\g\g',z)=h(\g,\g'z)\,h(\g',z)\qquad(\g,\g'\in\Gz{4\Ob}).
	\end{equation*}
	For $\b\in F$ we write $\ub{\beta}=\begin{psmallmatrix}1&\beta\\0&1\end{psmallmatrix}$; if $\b\in2\Diff^{-1}$ then $\theta(z+\beta)=\theta(z)$, so that
	\begin{equation}\label{eq:h-translation}
		h(\ub{\b},z)=1,\qquad h(\g\,\ub{\b},z)=h(\g,z+\b)
		\qquad(\b\in2\Diff^{-1},\ \g\in\Gz{4\Ob}).
	\end{equation}
	
	Let $k$ be an odd positive integer and put $(f|_k[\g])(z)=f(\g z)h(\g,z)^{-k}$. The weight is $k/2$.
	\begin{definition}\label{def:Mrho}
		Let $\mf$ be an integral ideal with $4\Ob\mid\mf$ and let $\rho:\Gz{\mf}\to\mathbb C^{\times}$ be a character of finite order. We write $M_k(\Gz{\mf};\rho)$ for the space of holomorphic functions $f:\HCal^{\ba}\to\mathbb C$ such that
		\begin{equation*}
			f|_k[\g]=\rho(\g)\,f\qquad(\g\in\Gz{\mf}).
		\end{equation*}
	\end{definition}

For a character $\psi$ of $(\Ob/\mf)^{\times}$ with $\psi(-1)=1$, we write $\Mk{\mf}{\psi}$ for the space corresponding to $\rho(\g)=\psi(d_\g)$. Every $f\in M_k(\Gz{\mf};\rho)$ such that $\rho(\ub{\b})=1$ for all $\b\in2\Diff^{-1}$ is invariant under translations by $2\Diff^{-1}$, by \eqref{eq:h-translation}, and therefore has a Fourier expansion
\begin{equation}\label{eq:fourier}
	f(z)=\sum_{\x\in2^{-1}\Ob_{+}\cup\{0\}}a_\x\,e(\x z),
\end{equation}
the indices being totally non-negative by the Koecher principle, since $g\ge2$. All the characters $\rho$ occurring below are trivial on the elements $\ub{\b}$, $\b\in2\Diff^{-1}$.
	
	\subsection{The Eisenstein series}\label{ss:eisenstein}
	For a totally imaginary quadratic extension $K/F$ let $\chi_K$ be the quadratic character of $K/F$ on ideals of $F$: $\chi_K(\pf)=1,-1,0$ according as the prime $\pf$ splits, is inert, or ramifies in $K$. Put $h^-(K/F)=h(K)/h(F)$. Let $w_K$ be the number of roots of unity in $K$ and $Q_K\in\{1,2\}$ its Hasse unit index. For $\x\in2^{-1}\Ob_{+}$ put
	\begin{equation*}
		\Kx{\x}=F(\sqrt{-2\x}),
		\qquad
		L_\x=\frac{2^{g}\,h^{-}(\Kx{\x}/F)}{Q_{\Kx{\x}}\,w_{\Kx{\x}}}.
	\end{equation*}
	We normalize $v_p$ by $v_p(p)=1$ and put $v_p(0)=+\infty$. For every odd prime $p$, we have $v_p(h^-(K/F))\ge0$. Put
	\begin{equation}\label{eq:beta}
		\etwo(\x)=\prod_{\qf\mid2}\bigl(1-\chi_{\Kx{\x}}(\qf)\bigr),
		\qquad
		\beta(\x)=\sum_{(\fa,\fb)}\mu(\fa)\,\chi_{\Kx{\x}}(\fa)\,\Norm(\fb),
	\end{equation}
	where $(\fa,\fb)$ runs over the pairs of integral ideals prime to $2\Ob$ with $(\fa\fb)^{2}\mid2\x\Ob$, and $\mu$ is the M\"obius function. The sum $\beta(\x)$ is a product over the primes $\pf\nmid2$ with $\pf^{2}\mid2\x\Ob$: if $\pf^{a}$ is the exact power of $\pf$ dividing $2\x\Ob$ and $n=\lfloor a/2\rfloor$, the factor at $\pf$ is
	\begin{equation}\label{eq:beta-local}
		\sum_{j=0}^{n}\Norm\pf^{\,j}-\chi_{\Kx{\x}}(\pf)\sum_{j=0}^{n-1}\Norm\pf^{\,j}
		\ \geq\ \Norm\pf^{\,n}.
	\end{equation}
	We put $b_E(\x)=\etwo(\x)\beta(\x)$.
	\begin{lemma}\label{lem:eisenstein}
		There is a form $\Es$ in $M_3(\Gz{4\Ob},1)$, of weight $3/2$, whose Fourier expansion \eqref{eq:fourier} has the coefficients
		\begin{equation*}
			a_\x(\Es)=b_E(\x)\,L_\x
			\qquad(\x\in2^{-1}\Ob_{+}).
		\end{equation*}
		The constant coefficient is $a_0(\Es)=\z_F(-1)\prod_{\qf\mid2}(1-\Norm\qf)$.
	\end{lemma}
	
\begin{proof}
	In \cite[Proposition~6.3]{Sh}, take
	\[
	m=1,\qquad k=3,\qquad \rho=\mu=0,\qquad
	\psi=1,\qquad \mathfrak c=4\Ob.
	\]
	Then $s_k=-1/2$ and $C=2^{g/2}\pi^g$. Put
	\[
	\Es(z)=C^{-1}
	\left.
	\left[
	L_{4\Ob}(4s+2,1)\,
	E'(z,s;3/2,1,4\Ob)
	\right]\right|_{s=-1/2}.
	\]
	The transformation law in \cite[\S4]{Sh}, together with
	\cite[Theorem~6.1 and Proposition~6.3]{Sh}, gives
	$\Es\in M_3(\Gz{4\Ob},1)$, with Fourier expansion
	\[
	\Es(z)=L_{4\Ob}(-1,1)
	+\sum_{\x\in2^{-1}\Ob_+}
	\beta(\x)L_{4\Ob}(0,\chi_{\Kx{\x}})\,e(\x z).
	\]
	
	Since $L_{4\Ob}$ omits the Euler factors above $2$, the
	relative class number formula gives
	\[
	L_{4\Ob}(0,\chi_{\Kx{\x}})
	=\etwo(\x)L(0,\chi_{\Kx{\x}})
	=\etwo(\x)L_\x.
	\]
	Thus $a_\x(\Es)=b_E(\x)L_\x$. The constant coefficient follows from
	\[
	L_{4\Ob}(-1,1)
	=\z_F(-1)\prod_{\qf\mid2}(1-\Norm\qf).
	\]
\end{proof}
	We also need coefficients that can be nonzero when a prime above $2$ splits. Say that $\alpha\in\Ob$ is a square modulo $4$ if $\alpha\equiv y^2\pmod{4\Ob}$ for some $y\in\Ob$. If $\eta\in\Ob_+$ and $-\eta\equiv y^2\pmod4$, then $(y+\sqrt{-\eta})/2$ generates an order of $K=F(\sqrt{-\eta})$ with discriminant $-\eta$. Hence
	\[
	\eta\Ob=\mathfrak F_\eta^2 D(K/F)
	\]
	for an integral ideal $\mathfrak F_\eta$. Define
	\begin{equation}\label{eq:cohen}
		\mathfrak C(\eta)=\sum_{\fa\mid\mathfrak F_\eta}
		\mu(\fa)\chi_K(\fa)\sigma_1(\mathfrak F_\eta\fa^{-1}),
		\qquad \sigma_1(\fb)=\sum_{\mathfrak c\mid\fb}\Norm\mathfrak c,
	\end{equation}
	and put
	\[
	b_G(\x)=
	\begin{cases}
		\mathfrak C(2\x),&-2\x\text{ is a square modulo }4,\\
		0,&\text{otherwise}.
	\end{cases}
	\]
	For every prime power $\pf^n$ exactly dividing $\mathfrak F_\eta$, the local factor of $\mathfrak C(\eta)$ is \eqref{eq:beta-local}, including at primes above $2$. In particular $b_G(\eta/2)=1$ if $\mathfrak F_\eta=\Ob$.

	\begin{lemma}\label{lem:plus-eisenstein}
		There is a form $\Gs\in M_3(\Gz{4\Ob},1)$, of weight $3/2$, with
		\[
		a_\x(\Gs)=b_G(\x)L_\x\quad(\x\in2^{-1}\Ob_+),\qquad
		a_0(\Gs)=\z_F(-1).
		\]
	\end{lemma}
	\begin{proof}
		Since $g\ge2$, \cite[Theorem~10.3]{Su} gives a holomorphic
		Eisenstein series $H(z)=G_{3/2}(z,1)$ with constant coefficient
		$\z_F(-1)$. Put $\Gs(z)=H(z/2)$.
		Su's theta series is $\theta(2z)$, and his group has
		$b\in\Diff^{-1}$ and $c\in4\Diff$. The change of variables
		$z\mapsto z/2$ gives our group $\Gz{4\Ob}$ with factor of
		automorphy $h^3$, and changes the Fourier index $\eta$ to $\eta/2$.
		The coefficient formula in that theorem, together with
		\eqref{eq:cohen} and $L(0,\chi_{\Kx{\x}})=L_\x$, gives
		$a_\x(\Gs)=b_G(\x)L_\x$. The constant coefficient is unchanged.
	\end{proof}
	
	\emph{Throughout we choose $\mathcal E=\Es$ or $\Gs$, put $b=b_E$ or $b_G$ accordingly, and write $a_\x=b(\x)L_\x$ for its positive-index Fourier coefficients.}
	
	The choice matters at $2$: $\Es$ vanishes on every field in which some prime above $2$ splits. Every CM quadratic field has an index in the support of $\Gs$: for $F(\sqrt{-\eta})$, take $2\x=4\eta$. This concerns nonvanishing; it does not assert that the coefficient is a unit at $p$ or that the new index satisfies a fixed local condition.
	
	\subsection{Local conditions}\label{ss:sieve}
	Let $\nf\neq\Ob$ be an integral ideal. Put $W_{\nf}=2\Diff^{-1}\nf^{-1}$, call $w\in W_{\nf}$ primitive if $w\notin2\Diff^{-1}\nf^{-1}\pf$ for every prime $\pf\mid\nf$, and put $\ew(x)=e(xw/2)$ for $x\in\Ob$. By \nAdd, if $w$ is primitive then $\ew$ is an additive character of $\Ob/\nf$ which is non-trivial on every ideal strictly larger than $\nf$. Consequently, for primitive $w$ and $c\in\Ob$,
	\begin{equation}\label{eq:orthogonality}
		\sum_{s\in\Ob/\nf}\ew(sc)=
		\begin{cases}
			\Norm\nf,&\text{if }c\in\nf,\\
			0,&\text{otherwise,}
		\end{cases}
	\end{equation}
	because $s\mapsto\ew(sc)$ is a character of $\Ob/\nf$, trivial if and only if $c\in\nf$.	
	\begin{proposition}[{\nPhase}]\label{prop:from-n}
		Let $\mf\subseteq4\Ob$ and $\nf\neq\Ob$ be integral ideals, let $\g=\begin{psmallmatrix}a&b\\c&d\end{psmallmatrix}\in\Gz{\mf\nf^{2}}$ and $w\in W_{\nf}$, and let $s,s'\in\Ob$ satisfy $s'\equiv d^{2}s\pmod{\nf}$. Then $\g_s:=\ub{sw}\,\g\,\ub{s'w}^{-1}$ lies in $\Gz{\mf\nf^{2}}$, its lower right entry is $d-s'wc\equiv d\pmod{\mf\nf}$, and
		\begin{equation*}
			h(\g,z)=h(\g_s,z+s'w)\qquad(z\in\HCal^{\ba}).
		\end{equation*}
	\end{proposition}

	The congruence for the lower right entry holds because $wc\in2\Diff^{-1}\nf^{-1}\cdot2^{-1}\mf\nf^{2}\Diff=\mf\nf$.
	
	We shall say that a subset $R\subseteq\Ob/\nf$ is \emph{stable under squares} if $u^{2}R=R$ for every unit $u$ of the ring $\Ob/\nf$. For $\eta\in\Ob$ we write $\eta\in R$ if the class of $\eta$ modulo $\nf$ belongs to $R$.		
	\begin{lemma}\label{lem:projection}
		Let $\mf\subseteq4\Ob$ and $\nf\neq\Ob$ be integral ideals, let $\psi$ be a character of $(\Ob/\mf)^{\times}$ with $\psi(-1)=1$, and let $R\subseteq\Ob/\nf$ be stable under squares. If $f\in\Mk{\mf}{\psi}$ has the expansion \eqref{eq:fourier}, then
		\begin{equation*}
			f_R(z)=\sum_{\x:\ 2\x\in R}a_\x\,e(\x z)
		\end{equation*}
		belongs to $\Mk{\mf\nf^{2}}{\psi}$.
	\end{lemma}
	
	\begin{proof}
		Choose a primitive $w\in W_{\nf}$. For $s\in\Ob$ we have $f(z+sw)=\sum_\x a_\x\,\ew(s\cdot2\x)\,e(\x z)$, and this function depends only on $s$ modulo $\nf$, because $\nf w\subseteq2\Diff^{-1}$. For $\alpha\in\Ob$ put
		\begin{equation}\label{eq:fa}
			f_\alpha(z)=\frac{1}{\Norm\nf}\sum_{s\in\Ob/\nf}\ew(-s\alpha)\,f(z+sw)
			=\sum_{\x:\ 2\x\equiv\alpha\ (\mathrm{mod}\ \nf)}a_\x\,e(\x z),
		\end{equation}
		where the second equality follows from \eqref{eq:orthogonality}. Thus $f_\alpha$ depends only on $\alpha$ modulo $\nf$, it is holomorphic, and $f_R=\sum_{\alpha\in R}f_\alpha$.
		
		Let $\g=\begin{psmallmatrix}a&b\\c&d\end{psmallmatrix}\in\Gz{\mf\nf^{2}}$. Since $ad-bc=1$ and $bc\in\mf\nf^{2}$, the element $d$ is a unit modulo $\nf$. For $s\in\Ob$ choose $s'\in\Ob$ with $s'\equiv d^{2}s\pmod{\nf}$. By Proposition~\ref{prop:from-n},
		\begin{equation*}
			f(\g z+sw)=f\bigl(\g_s(z+s'w)\bigr)
			=\psi(d)\,h(\g,z)^{k}\,f(z+s'w).
		\end{equation*}
		As $s$ runs through $\Ob/\nf$ so does $s'$, and $\ew(-s\alpha)=\ew(-s'\bar d^{\,2}\alpha)$, where $\bar d\in\Ob$ is an inverse of $d$ modulo $\nf$. Hence, by \eqref{eq:fa},
		\begin{equation*}
			f_\alpha(\g z)=\psi(d)\,h(\g,z)^{k}\,f_{\bar d^{2}\alpha}(z).
		\end{equation*}
		Summing over $\alpha\in R$ and using $\bar d^{\,2}R=R$ we obtain $f_R|_k[\g]=\psi(d)f_R$.
	\end{proof}
	\begin{definition}\label{def:local}
		A \emph{local condition} is a pair $(\nf,R)$ consisting of an integral ideal $\nf\neq\Ob$ and a non-empty subset $R\subseteq\Ob/\nf$ which is stable under squares and does not contain $0$.
	\end{definition}
	
	For an integral ideal $\mf$ with $4\Ob\mid\mf$ we put
	\begin{equation}\label{eq:kappa0}
		\Mo(\mf)\,\ZZ=\mf\Diff\cap\ZZ\quad(\Mo(\mf)>0),
		\qquad
		\ko(\mf)=\frac g8\,\idx(\Mo(\mf)),
	\end{equation}
	with $\idx(M)=M\prod_{u\mid M}(1+1/u)$, the product being over the primes $u$ dividing $M$, the index of $\Gamma_0(M)$ in $\SL_2(\ZZ)$. For a local condition $(\nf,R)$ we put
	\begin{equation}\label{eq:level}
		\mf=4\nf^{2},
		\qquad
		\Mo=\Mo(\mf),
		\qquad
		\ko=\ko(\mf).
	\end{equation}
	\begin{lemma}\label{lem:sieve}
		Let $\mathcal E$ be either series $\Es$ or $\Gs$ of Lemmas~\ref{lem:eisenstein} and~\ref{lem:plus-eisenstein}, with coefficients $a_\x$, and let $(\nf,R)$ be a local condition. The series
		\begin{equation*}
			\mathcal E_{R}(z)=\sum_{\x\in2^{-1}\Ob_{+},\ 2\x\in R}a_\x\,e(\x z)
		\end{equation*}
		belongs to $\Mk{\mf}{1}$ with $k=3$ and $\mf=4\nf^{2}$, and its constant term is $0$.
	\end{lemma}
	
	\begin{proof}
		Apply Lemma~\ref{lem:projection} to $f=\mathcal E$; the constant term disappears because $0\notin R$.
	\end{proof}
	
	We write $\aR_\x$ for the Fourier coefficients of $\mathcal E_R$: thus $\aR_\x=a_\x$ if $\x\in2^{-1}\Ob_{+}$ and $2\x\in R$, and $\aR_\x=0$ otherwise; in particular $\aR_0=0$.
	
	\begin{lemma}\label{lem:integrality}
		Let $\x\in2^{-1}\Ob_+$ and let $(\nf,R)$ be a local condition.
		\begin{enumerate}[label=\textup{(\alph*)},leftmargin=*]
			\item The number $b(\x)$ is a nonnegative integer and $a_\x\ge0$. We have $b_E(\x)\ne0$ if and only if no prime above $2$ splits in $K_\x$, and $b_G(\x)\ne0$ if and only if $-2\x$ is a square modulo $4$.
			\item If $p$ is odd and $K_\x\ne F(\z_p)$, then $a_\x\in\ZZ_{(p)}$. For a nonzero coefficient,
			\[
			0\le v_p(h^-(K_\x/F))\le v_p(a_\x).
			\]
			In particular $a_\x\not\equiv0\pmod p$ implies $p\nmid h^-(K_\x/F)$. The valuation inequality also holds whenever $a_\x\ne0$ and $K_\x/F$ ramifies at a prime not above $p$. There is an integer $c>0$ such that every coefficient of $c\mathcal E$ is integral.
			\item If $\lf=\pi\Ob$ is a prime ideal, $\pi\in\Ob_+$, and $\lf\nmid2$ and $\lf\nmid2\x$, then
			\[
			b(\pi^2\x)=(\Norm\lf+1-\chi_{K_\x}(\lf))b(\x).
			\]
			If also $\lf\nmid\nf$, then
			\[
			\aR_{\pi^2\x}=(\Norm\lf+1-\chi_{K_\x}(\lf))\aR_\x.
			\]
			\item If $\tau\in\Ob_+$ and $\tau\equiv1\pmod{8\Ob}$, then $b(\tau\x)\ne0$ if and only if $b(\x)\ne0$.
		\end{enumerate}
	\end{lemma}
	\begin{proof}
		For $\Es$, part (a) follows from \eqref{eq:beta-local} and
		$L_\x>0$. For $\Gs$, it follows directly from the definition
		of $b_G$ and $L_\x>0$.
		
		For (b), a quadratic CM extension contains a primitive $p$th root of unity exactly when it is $F(\z_p)$. Otherwise $p\nmid Q_Kw_K$, so $2^g/(Q_Kw_K)$ is a $p$-adic unit. The coefficient formula and $v_p(h^-(K/F))\ge0$ give the assertion. The ramification assertion follows because $F(\z_p)/F$ is unramified away from primes above $p$. There are finitely many possible $w_K$, since $\phi(w_K)\le2g$, where $\phi$ is Euler's totient function. Since $Q_K$ divides $2^g$, the denominator of $a_\x$ divides $w_K$. Thus the least common multiple of the possible $w_K$ clears every positive coefficient. Enlarge this integer to clear the rational constant term as well.
		
		For (c), $K_{\pi^2\x}=K_\x$, so $e_2$ is unchanged and $L_{\pi^2\x}=L_\x$. At an odd prime $\lf\nmid2\x$, multiplication by $\pi^2$ adds exactly the factor $1+\Norm\lf-\chi_{K_\x}(\lf)$ to either divisor sum. It also preserves the condition of being a square modulo $4$, since $\pi$ is a unit modulo $4$. If $\lf\nmid\nf$, then $\pi$ is a unit modulo $\nf$ and, since $R$ is stable under squares, $2\pi^2\x\in R$ if and only if $2\x\in R$. The relation for $\aR_\x$ follows.
		
		For (d), Hensel's lemma applied to $X^2-\tau$ at $1$ shows that $\tau$ is a square in every completion above $2$: $v_\qf(1-\tau)\ge3v_\qf(2)>2v_\qf(2)$. Hence the dyadic splitting conditions for $K_{\tau\x}$ and $K_\x$ are the same. Also $-2\tau\x\equiv-2\x\pmod4$, so the square-modulo-$4$ condition is unchanged.
	\end{proof}
	\section{The operators \texorpdfstring{$U_\pi$}{U} and \texorpdfstring{$V_\pi$}{V}}\label{sec:operators}
	
	Throughout this section $\pi\in\Ob_{+}$ is a totally positive element such that $\lf=\pi\Ob$ is a prime ideal, and $\mf$ is an integral ideal with $4\Ob\mid\mf$.
	
	\subsection{The character \texorpdfstring{$\kp$}{kappa}}
	For $\g=\begin{psmallmatrix}a&b\\c&d\end{psmallmatrix}\in\Gz{4\lf}$ put
	\begin{equation*}
		\g_\pi=\begin{pmatrix}\pi&0\\0&1\end{pmatrix}\g
		\begin{pmatrix}\pi&0\\0&1\end{pmatrix}^{-1}
		=\begin{pmatrix}a&\pi b\\ c/\pi&d\end{pmatrix}.
	\end{equation*}
	Since $\pi b\in2\Diff^{-1}$ and $c/\pi\in2^{-1}\cdot4\Diff$, we have $\g_\pi\in\Gz{4\Ob}$; more generally $\g_\pi\in\Gz{\mf}$ if $\g\in\Gz{\mf\lf}$. The map $\g\mapsto\g_\pi$ is a homomorphism, and
	\begin{equation}\label{eq:pi-conj}
		\pi\cdot\g z=\g_\pi(\pi z),
		\qquad
		j(\g_\pi,\pi z)=j(\g,z)
		\qquad(z\in\HCal^{\ba}),
	\end{equation}
	where $\pi z=(\pi_vz_v)_{v\in\ba}\in\HCal^{\ba}$ because $\pi$ is totally positive.
	
	\begin{lemma}\label{lem:kappa}
		For $\g\in\Gz{4\lf}$ the function
		\begin{equation*}
			\kp(\g)=\frac{h(\g_\pi,\pi z)}{h(\g,z)}
		\end{equation*}
		is independent of $z$, takes values in $\{\pm1\}$, and $\kp:\Gz{4\lf}\to\{\pm1\}$ is a homomorphism with $\kp(\ub{\b})=1$ for $\b\in2\Diff^{-1}$.
	\end{lemma}
	
	\begin{proof}
		Since $d_{\g_\pi}=d_\g$ and $j(\g_\pi,\pi z)=j(\g,z)$, the squares of $h(\g_\pi,\pi z)$ and $h(\g,z)$ coincide by \eqref{eq:h2}. Hence the quotient is a holomorphic function on $\HCal^{\ba}$ with values in $\{\pm1\}$, that is, a constant $\kp(\g)\in\{\pm1\}$. The cocycle relation for $h$, applied to $\g\g'$ and to $(\g\g')_\pi=\g_\pi\g'_\pi$, together with $\pi\cdot\g'z=\g'_\pi(\pi z)$, gives
		\begin{equation*}
			\kp(\g\g')
			=\frac{h(\g_\pi,\pi\cdot\g'z)}{h(\g,\g'z)}\cdot
			\frac{h(\g'_\pi,\pi z)}{h(\g',z)}
			=\kp(\g)\kp(\g').
		\end{equation*}
		If $\b\in2\Diff^{-1}$, then $\ub{\b}_\pi=\ub{\pi\b}$ and both values of $h$ equal $1$ by \eqref{eq:h-translation}.
	\end{proof}
	
	\subsection{The operator \texorpdfstring{$V_\pi$}{V}}
	For $f$ as in \eqref{eq:fourier} put
	\begin{equation*}
		(V_\pi f)(z)=f(\pi z)=\sum_{\x}a_\x\,e(\pi\x z).
	\end{equation*}
	
	\begin{lemma}\label{lem:V}
		Let $\psi$ be a character of $(\Ob/\mf)^{\times}$ with $\psi(-1)=1$. If $f\in\Mk{\mf}{\psi}$, then
		\begin{equation*}
			(V_\pi f)|_k[\g]=\psi(d_\g)\,\kp(\g)^{k}\;V_\pi f
			\qquad(\g\in\Gz{\mf\lf}).
		\end{equation*}
	\end{lemma}
	
	\begin{proof}
		Let $\g\in\Gz{\mf\lf}$. Then $\g_\pi\in\Gz{\mf}$, $d_{\g_\pi}=d_\g$, and by \eqref{eq:pi-conj} and Lemma~\ref{lem:kappa},
		\begin{equation*}
			f(\pi\cdot\g z)=f\bigl(\g_\pi(\pi z)\bigr)
			=\psi(d_\g)\,h(\g_\pi,\pi z)^{k}f(\pi z)
			=\psi(d_\g)\,\kp(\g)^{k}\,h(\g,z)^{k}\,(V_\pi f)(z).\qedhere
		\end{equation*}
	\end{proof}
	
	\subsection{The operator \texorpdfstring{$U_\pi$}{U}}
	For $f$ as in \eqref{eq:fourier} put
	\begin{equation*}
		(U_\pi f)(z)=\sum_{\x\in2^{-1}\Ob_{+}\cup\{0\}}a_{\pi\x}\,e(\x z).
	\end{equation*}
	Since $2\x\in\lf$ if and only if $\x=\pi\x'$ with $\x'\in2^{-1}\Ob$, we have
	\begin{equation}\label{eq:U-via-fp}
		(U_\pi f)(z)=f^{(\lf)}(z/\pi),
		\qquad
		f^{(\lf)}(z)=\sum_{\x:\ 2\x\in\lf}a_\x\,e(\x z).
	\end{equation}
	In the notation of Lemma~\ref{lem:projection}, $f^{(\lf)}=f_R$ for $\nf=\lf$ and $R=\{0\}$; hence $f^{(\lf)}\in\Mk{\mf\lf^{2}}{\psi}$ if $f\in\Mk{\mf}{\psi}$.
	
	\begin{lemma}\label{lem:U-theta}
		$U_\pi\theta=V_\pi\theta$.
	\end{lemma}
	
	\begin{proof}
		The coefficient of $e(\x z)$ in $\theta$ is the number of $x\in\Ob$ with $x^{2}=2\x$. If $x^{2}=2\pi\x$, then $x\in\lf$ because $\lf$ is prime; writing $x=\pi y$ gives $\pi y^{2}=2\x$. Hence $(U_\pi\theta)(z)
		=\sum_{\xi}\#\{y\in\Ob:\pi y^2=2\xi\}\,e(\xi z)
		=\sum_{y\in\Ob}e(\pi y^2z/2)
		=\theta(\pi z)=(V_\pi\theta)(z)$.
	\end{proof}
	
	\begin{proposition}\label{prop:U}
		Let $\psi$ be a character of $(\Ob/\mf)^{\times}$ with $\psi(-1)=1$, and let $f\in\Mk{\mf}{\psi}$. Then
		\begin{equation*}
			(U_\pi f)|_k[\g]=\psi(d_\g)\,\kp(\g)^{k}\;U_\pi f
			\qquad(\g\in\Gz{\mf\lf}).
		\end{equation*}
		Thus $U_\pi f$ and $V_\pi f$ belong to the same space $M_k(\Gz{\mf\lf};\rho)$, where $\rho(\g)=\psi(d_\g)\kp(\g)^{k}$.
	\end{proposition}
	
	\begin{proof}
		Let $\g=\begin{psmallmatrix}a&b\\c&d\end{psmallmatrix}\in\Gz{\mf\lf}$ with $b\in\pi\cdot2\Diff^{-1}$, and put $\g^{\pi}=\begin{psmallmatrix}a&b/\pi\\ \pi c&d\end{psmallmatrix}$. Then $\g^{\pi}\in\Gz{\mf\lf^{2}}$, $(\g^{\pi})_\pi=\g$, and $(\g z)/\pi=\g^{\pi}(z/\pi)$. Since $f^{(\lf)}\in\Mk{\mf\lf^{2}}{\psi}$, \eqref{eq:U-via-fp} gives
		\begin{equation*}
			(U_\pi f)(\g z)=f^{(\lf)}\bigl(\g^{\pi}(z/\pi)\bigr)
			=\psi(d)\,h(\g^{\pi},z/\pi)^{k}\,(U_\pi f)(z).
		\end{equation*}
		By the definition of $\kp$, applied to $\g^{\pi}\in\Gz{4\lf}$ at the point $z/\pi$, we have $h(\g,z)=\kp(\g^{\pi})\,h(\g^{\pi},z/\pi)$. Hence
		\begin{equation}\label{eq:U-step1}
			(U_\pi f)|_k[\g]=\psi(d)\,\kp(\g^{\pi})^{k}\;U_\pi f .
		\end{equation}
		
		We claim that $\kp(\g^{\pi})=\kp(\g)$ for every $\g\in\Gz{4\lf}$ with $b_\g\in\pi\cdot2\Diff^{-1}$. Indeed, $\theta\in M_1(\Gz{4\Ob},1)$ by the definition of $h$, so \eqref{eq:U-step1} with $f=\theta$, $k=1$, $\mf=4\Ob$ gives $(U_\pi\theta)|_1[\g]=\kp(\g^{\pi})\,U_\pi\theta$, while Lemma~\ref{lem:V} gives $(V_\pi\theta)|_1[\g]=\kp(\g)\,V_\pi\theta$. Since $U_\pi\theta=V_\pi\theta\neq0$ by Lemma~\ref{lem:U-theta}, the claim follows.
		
		Let $\g=\begin{psmallmatrix}a&b\\c&d\end{psmallmatrix}\in\Gz{\mf\lf}$ be arbitrary, and put $\b=db\in2\Diff^{-1}$ and $\g'=\g\,\ub{-\b}=\begin{psmallmatrix}a&b-a\b\\c&d-c\b\end{psmallmatrix} \in\Gz{\mf\lf}$. Since $ad-bc=1$ and $bc\in2\Diff^{-1}\cdot2^{-1}\mf\lf\Diff=\mf\lf$, we have $b-a\b=-b\cdot bc\in\pi\cdot2\Diff^{-1}$ and $d_{\g'}=d-c\b\equiv d\pmod{\mf}$. Moreover $\kp(\g')=\kp(\g)$ and $h(\g',z+\b)=h(\g,z)$ by Lemma~\ref{lem:kappa} and \eqref{eq:h-translation}. As $U_\pi f$ is invariant under $z\mapsto z+\b$, applying \eqref{eq:U-step1} to $\g'$ and using $\kp((\g')^\pi)=\kp(\g')$ gives
		\begin{align*}
			(U_\pi f)(\g z)=(U_\pi f)\bigl(\g'(z+\b)\bigr)
			&=\psi(d)\,\kp(\g')^{k}\,h(\g',z+\b)^{k}\,(U_\pi f)(z+\b)\\
			&=\psi(d)\,\kp(\g)^{k}\,h(\g,z)^{k}\,(U_\pi f)(z).
		\end{align*}
		Holomorphy of $U_\pi f$ is clear from \eqref{eq:U-via-fp}.
	\end{proof}

	\section{Twisted diagonal restriction}\label{sec:restriction}
	
	\subsection{Restriction along a totally positive integer}
	Let $\vp\in\Ob_{+}$. For a function $f$ on $\HCal^{\ba}$ and $z\in\HCal$ put
	\begin{equation*}
		(\Dv f)(z)=f(2\vp z),\qquad 2\vp z:=(2\vp_v z)_{v\in\ba}\in\HCal^{\ba}.
	\end{equation*}
	If $f$ has the expansion \eqref{eq:fourier}, then, with $\tq=e^{2\pi iz}$,
	\begin{equation}\label{eq:restricted-expansion}
		(\Dv f)(z)=\sum_{n\geq0}c_n\,\tq^{\,n},
		\qquad
		c_n=\sum_{\substack{\x\in2^{-1}\Ob_{+}\cup\{0\}\\ 2\Tr(\x\vp)=n}}a_\x .
	\end{equation}
	The inner sum is finite because $\vp$ is totally positive, and $2\Tr(\x\vp)\in\ZZ$ because $2\x\in\Ob$.
	
	\begin{lemma}\label{lem:restriction}
		Let $\mf$ be an integral ideal with $4\Ob\mid\mf$, let $\rho:\Gz{\mf}\to\mathbb C^{\times}$ be a character with $\rho^{4}=1$, and let $f\in M_k(\Gz{\mf};\rho)$. Let $\vp\in\Ob_{+}$ and let $M$ be a positive integer with $M\in\vp\,\mf\Diff$. Then $(\Dv f)^{4}$ is an elliptic modular form of weight $2kg$ on $\Gamma_0(M)\subseteq\SL_2(\ZZ)$ with trivial character.
	\end{lemma}
	
	\begin{proof}
		Let $\g=\begin{psmallmatrix}a&b\\c&d\end{psmallmatrix}\in\Gamma_0(M)$ and put
		\begin{equation*}
			\g_{\vp}=\begin{pmatrix}2\vp&0\\0&1\end{pmatrix}\g
			\begin{pmatrix}2\vp&0\\0&1\end{pmatrix}^{-1}
			=\begin{pmatrix}a&2\vp b\\ c/(2\vp)&d\end{pmatrix}.
		\end{equation*}
		Then $a,d\in\ZZ$, $2\vp b\in2\Ob\subseteq2\Diff^{-1}$, and $c/(2\vp)\in(2\vp)^{-1}M\ZZ\subseteq2^{-1}\mf\Diff$, so that $\g_{\vp}\in\Gz{\mf}$. Moreover $2\vp\cdot\g z=\g_{\vp}(2\vp z)$ and $j(\g_{\vp},2\vp z)=(cz+d)^{g}$. By \eqref{eq:h2} and $\rho^{4}=1$,
		\begin{equation*}
			(\Dv f)^{4}(\g z)=f\bigl(\g_{\vp}(2\vp z)\bigr)^{4}
			=\rho(\g_{\vp})^{4}\,j(\g_{\vp},2\vp z)^{2k}\,f(2\vp z)^{4}
			=(cz+d)^{2kg}\,(\Dv f)^{4}(z).
		\end{equation*}
	For every $\sigma\in\SL_2(\ZZ)$ we have
	\[
	(\Dv f)^4|_{2kg}\sigma
	=\Dv\bigl(f^4|_{2k}\sigma_{\vp}\bigr).
	\]
	By the Koecher principle, the Fourier indices $\nu$ of $f^4|_{2k}\sigma_{\vp}$ are totally non-negative. Since $\vp$ is totally positive, the restricted exponents $2\Tr(\nu\vp)$ are nonnegative, so $(\Dv f)^4$ is holomorphic at every cusp.
	\end{proof}
	
	For a power series $\phi=\sum_{n\geq0}c_n\tq^{\,n}$ with coefficients in $\ZZ_{(p)}$ we write $\ord_p(\phi)=\min\{n:\ c_n\not\equiv0\pmod p\}$, and $\ord_p(\phi)=\infty$ if $\phi\equiv0\pmod p$.
	
	\begin{corollary}\label{cor:sturm}
		Let $f$, $\vp$ and $M$ be as in Lemma~\ref{lem:restriction}, and suppose that $\rho(\ub{\b})=1$ for $\b\in2\Diff^{-1}$, so that $f$ has an expansion \eqref{eq:fourier}. Let $p$ be a prime, and suppose that there is a positive integer $c$ prime to $p$ such that $c\,a_\x\in\ZZ$ for all the coefficients $a_\x$ of $f$. If $\Dv f\not\equiv0\pmod p$, then
		\begin{equation*}
			\ord_p(\Dv f)\leq\frac{kg}{24}\,\idx(M).
		\end{equation*}
	\end{corollary}
	
	\begin{proof}
		The form $(c\,\Dv f)^{4}$ has integral coefficients, weight $2kg$ and level $M$, and it is not $\equiv0\pmod p$. By Sturm's theorem \cite{Sturm}, a form of weight $k'$ on $\Gamma_0(M)$ with integral coefficients which is not $\equiv0\pmod p$ has $\ord_p\leq\frac{k'}{12}\idx(M)$. Since $\FF_p[[\tq]]$ is an integral domain and $p\nmid c$, we have $\ord_p((c\,\Dv f)^{4})=4\ord_p(\Dv f)$, and the claim follows with $k'=2kg$.
	\end{proof}
	
	\subsection{The exponents for a prime of degree one}\label{ss:exponents}
	Now let $\lf$ be a prime ideal of $F$ of degree one which is unramified over $\QQ$, so that $\ell=\Norm\lf$ is a prime number and $\lf\nmid\Diff$, and suppose that $\lf=\pi\Ob$ with $\pi\in\Ob_{+}$. We put
	\begin{equation*}
		\vp=\ell/\pi.
	\end{equation*}
	
	\begin{lemma}\label{lem:exponents}
		\begin{enumerate}[label=\textup{(\alph*)},leftmargin=*]
			\item $\vp\in\Ob_{+}$, $\vp\,\lf=\ell\Ob$, and $\Norm_{F/\QQ}(\vp)=\ell^{g-1}$.
			\item For $\nu\in2^{-1}\Ob$ we have $2\Tr(\pi\nu\vp)=2\ell\Tr(\nu)$.
			\item For $\x\in2^{-1}\Ob$ we have $\ell\mid2\Tr(\x\vp)$ if and only if $2\x\in\lf$.
			\item If $\x\in2^{-1}\Ob_{+}$ and $2\Tr(\x\vp)=n$, then $\Norm_{F/\QQ}(2\x)\leq(n/g)^{g}\,\ell^{1-g}$.
		\end{enumerate}
	\end{lemma}
	
	\begin{proof}
		(a) Since $\ell=\Norm_{F/\QQ}(\pi)$, the quotient $\ell/\pi$ is an algebraic integer of $F$; it is totally positive because $\ell$ and $\pi$ are, and $\Norm(\vp)=\ell^{g}/\ell$. (b) is clear from $\pi\vp=\ell$.
		
		(c) The map $\phi:\Ob\to\ZZ/\ell\ZZ$, $\eta\mapsto\Tr(\eta\vp)\bmod\ell$, is additive and vanishes on $\lf=\pi\Ob$ by (b). Since $\Ob/\lf$ has prime order $\ell$, either $\ker\phi=\lf$ or $\phi=0$. If $\phi=0$, then $\Tr(\eta\,\vp/\ell)\in\ZZ$ for all $\eta\in\Ob$, that is, $\pi^{-1}=\vp/\ell\in\Diff^{-1}$, so $\Diff\subseteq\pi\Ob=\lf$, contrary to $\lf\nmid\Diff$. Hence $\ker\phi=\lf$; apply this to $\eta=2\x$.
		
		(d) Since $2\x\vp$ is totally positive with trace $n$, the AM--GM inequality gives $\Norm(2\x\vp)=\Norm(2\x)\Norm(\vp)\leq(n/g)^g.$ The claim follows from (a).
	\end{proof}
	
	By (b) and (c), the coefficients $c_n$ of \eqref{eq:restricted-expansion} are of two kinds: those with $\ell\nmid n$ only involve indices $\x$ with $\lf\nmid2\x$, and
	\begin{equation}\label{eq:old-exponents}
		c_{\ell m}=\sum_{\substack{\nu\in2^{-1}\Ob_{+}\cup\{0\}\\ 2\Tr(\nu)=m}}a_{\pi\nu}
		\qquad(m\geq0).
	\end{equation}
	
	\subsection{Production of a new index}\label{ss:production}
	Let $\lf=\pi\Ob$ and $\vp=\ell/\pi$ be as in \S\ref{ss:exponents}. Let $\mf$ be an integral ideal with $4\Ob\mid\mf$ such that $\ell\nmid\Mo(\mf)$; since $\Mo(\mf)\in\mf\Diff$, this implies $\lf\nmid\mf\Diff$, in particular $\lf\nmid2$. Let $\psi$ be a character of $(\Ob/\mf)^{\times}$ with $\psi(-1)=1$ and $\psi^{4}=1$, let $f\in M_3(\Gz{\mf},\psi)$, with rational Fourier coefficients $a_\x$ as in \eqref{eq:fourier}, and let $p$ be a prime. For $\alpha\in\ZZ_{(p)}$ put
	\begin{equation}\label{eq:H}
		H_\pi=U_\pi f-\alpha\,V_\pi f,
		\qquad
		\Phi=\Dv H_\pi=\sum_{n\geq0}c_n\,\tq^{\,n}.
	\end{equation}
	By Lemma~\ref{lem:V} and Proposition~\ref{prop:U}, $H_\pi\in M_3(\Gz{\mf\lf};\rho)$ with $\rho(\g)=\psi(d_\g)\kp(\g)^{3}$, and its Fourier coefficients are
	\begin{equation}\label{eq:b}
		a_\x(H_\pi)=a_{\pi\x}-\alpha\,a_{\x/\pi},
	\end{equation}
	where $a_{\x/\pi}=0$ if $\x/\pi\notin2^{-1}\Ob$. If $\alpha=\alpha_1/\alpha_2$ with integers $\alpha_1$, $\alpha_2$ and $p\nmid\alpha_2$, and if a positive integer $c$ prime to $p$ clears all $a_\x$, then $c\alpha_2\,a_\x(H_\pi)\in\ZZ$ for all $\x$. The first two identities below require no integrality assumption.
	\begin{lemma}\label{lem:coefficients}
		\begin{enumerate}[label=\textup{(\alph*)},leftmargin=*]
			\item If $\ell\nmid n$, then $c_n=\sum a_{\pi\x}$, the sum being over the $\x\in2^{-1}\Ob_{+}$ with $2\Tr(\x\vp)=n$; all these $\x$ satisfy $\lf\nmid2\x$.
			\item For every $m\geq0$,
			\begin{equation*}
				c_{\ell m}=\sum_{\substack{\nu\in2^{-1}\Ob_{+}\cup\{0\}\\ 2\Tr(\nu)=m}}
				\bigl(a_{\pi^{2}\nu}-\alpha\,a_\nu\bigr).
			\end{equation*}
			\item Suppose that a positive integer $c$ prime to $p$ satisfies $c a_\x\in\ZZ$ for every index, including $\x=0$. If $\Phi\not\equiv0\pmod p$, then $\ord_p(\Phi)\leq\ko(\mf)(\ell+1)$.
		\end{enumerate}
	\end{lemma}
	
	\begin{proof}
		(a) By Lemma~\ref{lem:exponents}\textup{(c)}, $2\Tr(\x\vp)=n$ with $\ell\nmid n$ implies $2\x\notin\lf$, so $\x\neq0$, $\x/\pi\notin2^{-1}\Ob$, and $a_\x(H_\pi)=a_{\pi\x}$ by \eqref{eq:b}.
		
		(b) By \eqref{eq:old-exponents}, $c_{\ell m}=\sum a_{\pi\nu}(H_\pi)$ over the $\nu\in2^{-1}\Ob_{+}\cup\{0\}$ with $2\Tr(\nu)=m$, and $a_{\pi\nu}(H_\pi)=a_{\pi^{2}\nu}-\alpha\,a_\nu$ by \eqref{eq:b}.
		
		(c) We apply Corollary~\ref{cor:sturm} to $H_\pi$, with $\mf\lf$, $\rho$ as above, the integer $c\alpha_2$ and $M=\ell\Mo(\mf)$; indeed $\rho^{4}=1$ because $\psi^{4}=1$ and $\kp$ takes the values $\pm1$, $\rho(\ub{\b})=1$ for $\b\in2\Diff^{-1}$ by Lemma~\ref{lem:kappa}, and $\ell\Mo(\mf)\in\ell\,\mf\Diff=\vp\,\mf\lf\,\Diff$ by Lemma~\ref{lem:exponents}\textup{(a)}. As $\ell\nmid\Mo(\mf)$, $\frac{3g}{24}\,\idx(\ell\Mo(\mf))=\frac g8\,\idx(\Mo(\mf))(\ell+1)=\ko(\mf)(\ell+1)$.
	\end{proof}
	\begin{proposition}[Production with a valuation bound]
		\label{prop:valuation-production}\label{prop:production}
		Let $p$ be a prime, let $4\Ob\mid\mf$, and let $\psi$ be a character of $(\Ob/\mf)^\times$ with $\psi(-1)=1$ and $\psi^4=1$. Let $f\in M_3(\Gz\mf;\psi)$ have rational Fourier coefficients, and suppose that its Fourier coefficients $a_\x$, for $\x\in2^{-1}\Ob_+$, have a common denominator prime to $p$. Let $\alpha\in\ZZ_{(p)}$, and let $\lf=(\pi)$ have degree one, with $\pi\in\Ob_+$ and
		\[
		\ell=\Norm\lf,\qquad
		\ell\nmid\Mo(\mf),\qquad
		\ell>8\ko(\mf).
		\]
		Write
		\[
		D_{\ell/\pi}(U_\pi-\alpha V_\pi)f
		=\sum_{n\ge0}c_n\tq^n.
		\]
		Suppose that an integer $b\ge0$ satisfies
		\[
		v_p(c_{\ell j})\ge b+1
		\qquad(0\le j\le\ko(\mf),\ j\in\ZZ),
		\]
		and that
		\[
		v_p(c_{n_*})\le b
		\qquad\text{for some integer }n_*\ge1.
		\]
		Then there is $\x\in2^{-1}\Ob_+$ with
		\[
		\lf\nmid2\x,\qquad
		v_p(a_{\pi\x})\le b,\qquad
		\Norm(2\x)\le(2\ko(\mf)/g)^g\ell.
		\]
		The field $K=F(\sqrt{-2\pi\x})$ is ramified at $\lf$ and satisfies
		\[
		\Norm D(K/F)\le(8\ko(\mf)/g)^g\ell^2.
		\]
	\end{proposition}

	\begin{proof}
		Put $\kappa=\ko(\mf)$ and $\Phi=\sum_{n\ge0}c_nq^n$, where $q=e^{2\pi iz}$. Every positive-index coefficient of $\Phi$ is $p$-integral, and the hypothesis with $j=0$ gives the same property for its constant term. Thus
		\[
		h=\min_{n\ge0}v_p(c_n)
		\quad\text{satisfies}\quad
		0\le h\le v_p(c_{n_*})\le b.
		\]
		The series $p^{-h}\Phi$ has a common denominator prime to $p$ and nonzero reduction. Lemma~\ref{lem:restriction} shows that its fourth power is an elliptic modular form of weight $6g$ and level $\ell\Mo(\mf)$. Clear its denominator by an integer prime to $p$ and apply Sturm's theorem, exactly as in the proof of Corollary~\ref{cor:sturm}. The first coefficient of $\Phi$ of valuation $h$ occurs at an index
		\[
		n\le\frac{6g}{48}\idx(\ell\Mo(\mf))=\kappa(\ell+1).
		\]
		This argument requires integrality of the normalized restricted series, not of the normalized Hilbert form.
		
		If $n=\ell j$, then $j<\kappa+1/8$. Since $8\kappa$ and $8j$ are integers, $j\le\kappa$, contrary to $v_p(c_{\ell j})\ge b+1>h$. Hence $\ell\nmid n$. The coefficient identity in Lemma~\ref{lem:coefficients}(a) gives
		\[
		c_n=\sum_{2\Tr(\x\ell/\pi)=n}a_{\pi\x}.
		\]
		At least one summand has valuation at most $h$, and every index in this sum satisfies $\lf\nmid2\x$. By Lemma~\ref{lem:exponents}(d),
		\[
		\Norm(2\x)\le\left(\frac{\kappa(\ell+1)}g\right)^g\ell^{1-g}
		=\left(\frac\kappa g\right)^g(1+1/\ell)^g\ell
		\le\left(\frac{2\kappa}g\right)^g\ell.
		\]
		Finally $\ord_\lf(2\pi\x)=1$, so $K$ ramifies at $\lf$. The order $\Ob[\sqrt{-2\pi\x}]$ has discriminant $-8\pi\x$; therefore
		\[
		\Norm D(K/F)\le4^g\ell\Norm(2\x)
		\le(8\kappa/g)^g\ell^2.
		\]
	\end{proof}

	For $b=0$, the valuation hypotheses say that the small old coefficients vanish modulo $p$ and that some coefficient is nonzero modulo $p$. By Lemma~\ref{lem:coefficients}(b), it is enough to have $a_{\pi^2\nu}\equiv\alpha a_\nu\pmod p$ for $2\Tr\nu\le\ko(\mf)$, including $\nu=0$, and a nonzero trace sum of $a_{\pi^2\nu}-\alpha a_\nu$ modulo $p$ at some $m_0\ge1$. The produced coefficient is then a $p$-unit.
	\section{The main theorem}\label{sec:main}
	
	\subsection{Statement}
	Let $F$ be a totally real number field of degree $g\geq2$, choose $\mathcal E\in\{\Es,\Gs\}$, and let $(\nf,R)$ be a local condition. Let $\mf$, $\Mo$, $\ko$ and the coefficients $\aR_\x$ of $\mathcal E_R$ be as in \S\ref{ss:sieve}. We use the following condition.
	
	\begin{itemize}[leftmargin=3.2em]
		\item[\textup{(W)}] There exist an integer $m_0\ge1$ and a prime ideal $\Qf\nmid2$ with $\Norm\Qf>(\ko/g)^g$ such that
		\begin{equation}\label{eq:W}
			w:=\sum_{\substack{\x\in2^{-1}\Ob_+,\ 2\Tr(\x)=m_0\\
					\ord_\Qf(2\x)\ \mathrm{odd}}}\aR_\x\ne0.
		\end{equation}
	\end{itemize}
	
	Here nonvanishing is in $\QQ$, rather than modulo a specified prime.
	
	\begin{lemma}\label{lem:witness}
		Condition \textup{(W)} holds if and only if $\aR_\x\ne0$ for some $\x\in2^{-1}\Ob_+$. Equivalently, there is $\eta\in\Ob_+\cap R$ with $b(\eta/2)\ne0$.
	\end{lemma}
	
	\begin{proof}
		If \textup{(W)} holds, at least one summand in \eqref{eq:W} is nonzero. Conversely, suppose that $\eta\in\Ob_+\cap R$ and $b(\eta/2)\ne0$. Choose a prime ideal $\Qf=(\tau)\nmid2\eta$ with $\tau\in\Ob_+$, $\tau\equiv1\pmod{8\nf}$ and $\Norm\Qf>(\ko/g)^g$. Primes splitting in the narrow ray class field modulo $8\nf$ supply infinitely many choices. The element $\eta\tau$ belongs to $R$, has valuation one at $\Qf$, and satisfies $b(\eta\tau/2)\ne0$ by Lemma~\ref{lem:integrality}(d). Taking $m_0=\Tr(\eta\tau)$ puts the positive term $\aR_{\eta\tau/2}$ in \eqref{eq:W}; all other terms are nonnegative.
	\end{proof}
	
	For an odd prime $p$ and an integer $t\ge0$, let $N_{p,t}(X)$ count the distinct fields in
	\[
	\left\{K=F(\sqrt{-\eta})\ \middle|\
	\begin{aligned}
		&\eta\in\Ob_+,\quad\eta\in R,\\
		&\Norm D(K/F)<X,\quad v_p(h^-(K/F))\le t
	\end{aligned}\right\}.
	\]
	\begin{maintheorem}\hypertarget{main-theorem}{}
		Assume \textup{(W)}, and fix $m_0,\Qf,w$ as above. Then $w\in\ZZ[1/2]$. There are constants $C>0$ and $X_0>1$, depending on $F,\mathcal E,(\nf,R),m_0,\Qf$, but not on $p$, such that
		\[
		N_{p,v_p(w)}(X)\ge C\frac{\sqrt X}{\log X}
		\qquad\text{for every odd prime $p$ and every $X\ge X_0$.}
		\]
		In particular, if $v_p(w)=0$, then
		\[
		\#\left\{K=F(\sqrt{-\eta})\ \middle|\
		\begin{aligned}
			&\eta\in\Ob_+,\quad\eta\in R,\\
			&p\nmid h^-(K/F),\quad\Norm D(K/F)<X
		\end{aligned}\right\}
		\ge C\frac{\sqrt X}{\log X}\qquad(X\ge X_0).
		\]
	\end{maintheorem}
	For every odd prime outside the finite set dividing the numerator of $w$, the fields counted satisfy $p\nmid h^-(K/F)$. The fields produced may depend on $p$.
	
	The following example shows why we use both $\Es$ and $\Gs$: for a fixed local condition, one series may have nonempty retained support while the other has none.
	\begin{example}
		Let $F=\QQ(\sqrt2)$, $\Ob=\ZZ[\sqrt2]$, and put
		\[
		\nf=8\Ob,\qquad
		R_\pm=\{\pm u^2:u\in(\Ob/8\Ob)^\times\}.
		\]
		Both $(\nf,R_\pm)$ are local conditions. If $\eta\in\Ob_+\cap R_-$, then $-\eta$ is a square in the completion at the unique prime above $2$. Hence $\Es_{R_-}=0$. On the other hand, $7\in R_-$ and $-7\equiv1\pmod{4\Ob}$, so $a_{7/2}(\Gs)>0$. Conversely, $-1$ is not a square modulo $4\Ob$, so neither is $-\eta$ for any $\eta\in R_+$. Thus $\Gs_{R_+}=0$. Since the prime above $2$ does not split in $F(i)$, we have $a_{1/2}(\Es)>0$, and $1\in R_+$. Lemma~\ref{lem:witness} therefore supplies witnesses using $\Gs$ for $R_-$ and $\Es$ for $R_+$. In each case the other series has no retained positive coefficient.
	\end{example}
	The indivisibility assertion requires only that the witness be a $p$-adic unit. The individual retained coefficients need not all be $p$-integral.

	\subsection{Proof of the Main Theorem}
	\label{ss:proof-main}

	Let $m_0$ and $\Qf$ be as in \textup{(W)}, put $S=\{\nu\in2^{-1}\Ob_+:2\Tr\nu\le\ko,\ \aR_\nu\ne0\}$, and let $T=\{\nu\in2^{-1}\Ob_{+}:2\Tr(\nu)=m_0,\ \aR_\nu\neq0\}$. Let $\Delta\subseteq F^{\times}/(F^{\times})^{2}$ be the finite subgroup generated by the classes of $-2\nu$, $\nu\in S\cup T$, let $L_1=F(\sqrt{\delta}:\delta\in\Delta)$, let $H^{+}$ be the narrow Hilbert class field of $F$, and $L=H^{+}L_1$; this is a finite abelian extension of $F$.
	
	The map $\lambda(\delta)=(-1)^{\ord_{\Qf}(\delta)}$ is a well-defined character of $\Delta$. By Kummer theory there is a unique $\sigma_1\in\Gal(L_1/F)$ with $\sigma_1(\sqrt\delta)=\lambda(\delta)\sqrt\delta$ for all $\delta\in\Delta$. We claim that $\sigma_1$ is trivial on $L_1\cap H^{+}$. Indeed, $L_1\cap H^{+}=F(\sqrt\delta:\delta\in\Delta_0)$ with $\Delta_0=\{\delta\in\Delta:F(\sqrt\delta)\subseteq H^{+}\}$; for $\delta\in\Delta_0$ the extension $F(\sqrt\delta)/F$ is unramified at $\Qf$, which for $\Qf\nmid2$ means that $\ord_{\Qf}(\delta)$ is even, so $\lambda(\delta)=1$. Consequently there is $\sigma\in\Gal(L/F)$ which is trivial on $H^{+}$ and restricts to $\sigma_1$ on $L_1$.
	
	Put $B=(8\ko/g)^{g}$ and
	\begin{equation*}
		\ell_0=\max\bigl\{8\ko,\ B,\ (m_0/g)^{g}\bigr\},
	\end{equation*}
	and let $\mathcal P$ be the set of prime ideals $\lf$ of $F$ of degree one, unramified in $L$, with Frobenius $\sigma$ in $\Gal(L/F)$, and such that $\ell=\Norm\lf$ satisfies $\ell\nmid\Mo$ and $\ell>\ell_0$. By the Chebotarev density theorem, and since the number of prime ideals of $F$ of degree at least two and norm at most $Y$ is $O(Y^{1/2})$,
	\begin{equation}\label{eq:chebotarev}
		\#\{\lf\in\mathcal P:\ \Norm\lf<Y\}
		=\Bigl(\frac{1}{[L:F]}+o(1)\Bigr)\frac{Y}{\log Y}
		\qquad(Y\to\infty).
	\end{equation}
	
	Choose an integer $c>0$ clearing the coefficients of $\mathcal E_R$, as in Lemma~\ref{lem:integrality}(b). All the preceding choices are independent of $p$.
	
	\emph{Producing a field from an auxiliary prime.} Fix an odd prime $p$, and take $\lf\in\mathcal P$ with $\lf\nmid p$. Since $\lf$ splits completely in $H^+$, we have $\lf=\pi\Ob$ with $\pi\in\Ob_+$. For $\nu\in S$ we have $\Norm(2\nu)\leq(2\Tr(\nu)/g)^g\leq(\ko/g)^g<\Norm\Qf$, so $\ord_\Qf(2\nu)=0$, $\sigma$ fixes $\sqrt{-2\nu}$, and $\lf$ splits in $K_\nu$. For $\nu\in T$ we get $\chi_{K_\nu}(\lf)=(-1)^{\ord_\Qf(2\nu)}$.
	
	For positive $\nu$ with $2\Tr\nu\le\ko$ or $2\Tr\nu=m_0$, the size conditions give $\Norm(2\nu)<\ell$; here $B=(8\ko/g)^g>(\ko/g)^g$. Lemma~\ref{lem:integrality}(c) therefore gives
	\[
	\aR_{\pi^2\nu}-\ell\aR_\nu
	=(1-\chi_{K_\nu}(\lf))\aR_\nu.
	\]
	Write $D_{\ell/\pi}(U_\pi-\ell V_\pi)(c\mathcal E_R) =\sum_{n\ge0}c_n\tq^n$. Lemma~\ref{lem:coefficients}(b) now gives
	\[
	c_{\ell j}=0\quad(0\le j\le\ko,\ j\in\ZZ),\qquad
	c_{\ell m_0}=2cw.
	\]
	The constant coefficient is zero because $0\notin R$. Since $cw$ is a nonzero integer, we may apply Proposition~\ref{prop:production} with $b=v_p(cw)\ge0$. It gives a field $K(\lf)=F(\sqrt{-2\pi\x})$ ramified at $\lf$, with
	\[
	\Norm D(K(\lf)/F)\le B\ell^2,\qquad
	v_p(c\aR_{\pi\x})\le v_p(cw).
	\]
	The coefficient is nonzero, so $2\pi\x\in R$ and $b(\pi\x)\ne0$. Cancelling $v_p(c)$ and using $\lf\nmid p$ in Lemma~\ref{lem:integrality}(b), we obtain
	\[
	0\le v_p(h^-(K(\lf)/F))\le v_p(\aR_{\pi\x})\le v_p(w).
	\]
	There are infinitely many eligible auxiliary primes for each $p$. Consequently $v_p(w)\ge0$ for every odd $p$, and $w\in\ZZ[1/2]$. The field $K(\lf)$ is counted by $N_{p,v_p(w)}(X)$ whenever $B\ell^2<X$.
	
	It remains to bound the multiplicity of $\lf\mapsto K(\lf)$. Suppose that $K(\lf_1)=K(\lf_2)=K(\lf_3)=K$ for distinct $\lf_j\in\mathcal P$ of norms $\ell_1\leq\ell_2\leq\ell_3$. Each $\lf_j$ ramifies in $K$, so $\lf_1\lf_2\lf_3\mid D(K/F)$ and $\ell_1\ell_2\ell_3\leq\Norm(D(K/F))\leq B\ell_1^{2}$, whence $\ell_3\leq B$, contrary to $\ell_3>\ell_0\geq B$. Thus the map is at most two-to-one. Let $\pi_{\mathcal P}(Y)$ count primes in $\mathcal P$ of norm strictly less than $Y$. Excluding primes above $p$ removes at most $g$ degree-one primes. Hence
	\[
	N_{p,v_p(w)}(X)\ge
	\frac12\max\{0,\pi_{\mathcal P}(\sqrt{X/B})-g\}.
	\]
	By \eqref{eq:chebotarev}, this gives constants $C>0$ and $X_0>1$ independent of $p$. This proves the \hyperlink{main-theorem}{Main Theorem}. \qed
	
	\subsection{Choosing and combining witnesses}\label{ss:witnesses}
	
	In the following cases the sum \eqref{eq:W} reduces to a single term.
	
	\begin{corollary}\label{cor:q}
		\begin{enumerate}[label=\textup{(\alph*)},leftmargin=*]
			\item Let $\tau\in\Ob_+$ be such that $\Qf=\tau\Ob$ is a prime ideal not dividing $2$, with $\Norm\Qf>(\ko/g)^g$, and such that
			\begin{equation}\label{eq:balanced}
				(g+1)v(\tau)>\Tr(\tau)\qquad\text{for every real place $v$ of $F$.}
			\end{equation}
			Suppose that $\tau\in R$ and $b(\tau/2)\ne0$. Then \textup{(W)} holds with $m_0=\Tr(\tau)$ and $w=\aR_{\tau/2}$. In particular, this applies to $\tau=q$ for every odd rational prime $q>\ko/g$ which is inert in $F$ and satisfies $q\in R$ and $b(q/2)\ne0$.
			\item Let $g=2$, and let $q$ be an odd rational prime which splits in $F$, with $q>(\ko/2)^2$, $q\in R$ and $b(q/2)\ne0$. Then \textup{(W)} holds with $m_0=2q$, $\Qf$ a prime above $q$, and $w=\aR_{q/2}$.
		\end{enumerate}
		Put $K_*=F(\sqrt{-\tau})$ in \textup{(a)} and $K_*=F(\sqrt{-q})$ in \textup{(b)}. For every odd prime $p$, set $t_p=v_p(h^-(K_*/F))$. Then, with constants independent of $p$,
		\[
		N_{p,t_p}(X)\ge C\frac{\sqrt X}{\log X}\qquad(X\ge X_0).
		\]
		In particular, $p\nmid h^-(K_*/F)$ gives the indivisibility conclusion of the \hyperlink{main-theorem}{Main Theorem}.
	\end{corollary}
	
	\begin{proof}
		
		(a) Let $m_0=\Tr(\tau)$, and let $\x\in2^{-1}\Ob_{+}$ with $2\Tr(\x)=m_0$ and $\eta:=2\x\in\Qf$. Then $\eta=\tau\mu$ with $\mu\in\Ob_{+}$. We claim that $\mu=1$. By the inequality between the arithmetic and geometric means, $\Tr(\mu)\geq g\,\Norm(\mu)^{1/g}\geq g$, and $\Tr(\mu)=g$ only if all the conjugates of $\mu$ are equal to $1$. Hence, if $\mu\neq1$, the integer $\Tr(\mu)$ is at least $g+1$, and by \eqref{eq:balanced}
		\begin{equation*}
			\Tr(\eta)=\sum_{v}v(\tau)\,v(\mu)\geq\min_v v(\tau)\cdot\Tr(\mu)
			\geq(g+1)\min_v v(\tau)>\Tr(\tau),
		\end{equation*}
		a contradiction. Thus the sum \eqref{eq:W} is the single term $\aR_{\tau/2}$; note that $\ord_{\Qf}(\tau)=1$. If $q$ is a rational prime which is inert in $F$, then $\tau=q$ generates the prime ideal $\Qf=q\Ob$ of norm $q^{g}$, and \eqref{eq:balanced} holds because all the conjugates of $q$ are equal.
		
		(b) Let $q\Ob=\Qf\Qf'$, so that $\Norm\Qf=q>(\ko/2)^{2}$, and let $\eta\in\Ob_{+}\cap\Qf$ with $\Tr(\eta)=2q$. The conjugate $2q-\eta$ of $\eta$ lies in $\Qf$ and is the image of $\eta$ under the automorphism which exchanges $\Qf$ and $\Qf'$; hence $\eta\in\Qf\cap\Qf'=q\Ob$, and $\eta=q\mu$ with $\mu\in\Ob_{+}$ and $\Tr(\mu)=2$, so that $\mu=1$ as in \textup{(a)}. Thus the sum \eqref{eq:W} is the single term $\aR_{q/2}$, with $\ord_{\Qf}(q)=1$.
		
		In both cases the ideal generated by the index $\eta$ ($\eta=\tau$ or $q$) is squarefree and prime to $2$. For $\Es$, this gives $\beta(\eta/2)=1$, so $b_E(\eta/2)$ is a power of $2$ when nonzero. For $\Gs$, the identity $\eta\Ob=\mathfrak F_\eta^2D(K_*/F)$ gives $\mathfrak F_\eta=\Ob$, so $b_G(\eta/2)=1$. Thus the single coefficient $w$ is positive and, for every odd $p$,
		\[
		v_p(w)=v_p(h^-(K_*/F))-v_p(w_{K_*})\le t_p.
		\]
		The \hyperlink{main-theorem}{Main Theorem} gives $v_p(w)\ge0$ and the asserted count. This also covers $K_*=F(\z_p)$; no exclusion on the seed field is needed.
	\end{proof}
	
	\begin{definition}\label{def:witness-ideal}
		Fix $F,\mathcal E,(\nf,R)$ as above, and assume \textup{(W)}. Let $\mathcal W$ be the set of nonzero sums \eqref{eq:W}, as $m_0\ge1$ and $\Qf\nmid2$ vary with $\Norm\Qf>(\ko/g)^g$. Put
		\[
		I_w=I_w(F,\mathcal E;\nf,R)
		=(w:w\in\mathcal W)\subseteq\ZZ[1/2].
		\]
	\end{definition}
	
	The \hyperlink{main-theorem}{Main Theorem} shows that $I_w$ is a nonzero ideal of $\ZZ[1/2]$. Write $I_w=d_w\ZZ[1/2]$, where $d_w$ is a positive odd integer.
	
	\begin{corollary}\label{cor:witness-gcd}
		Assume \textup{(W)}. There are constants $C>0$ and $X_0>1$, independent of $p$, such that
		\[
		N_{p,v_p(d_w)}(X)\ge C\frac{\sqrt X}{\log X}
		\qquad\text{for every odd prime $p$ and every $X\ge X_0$.}
		\]
		In particular, $p\nmid d_w$ gives the indivisibility conclusion of the \hyperlink{main-theorem}{Main Theorem}.
		
		For any $w_1,\ldots,w_r\in\mathcal W$, let $D$ be the odd part of the greatest common divisor of their numerators in lowest terms. Then $d_w\mid D$, and the same assertion holds with $D$ in place of $d_w$, with constants depending on the chosen witnesses. If $D=1$, then $I_w=\ZZ[1/2]$ and the indivisibility conclusion holds for every odd prime.
	\end{corollary}
	
	\begin{proof}
		For a finite set of witnesses as above, we have
		\[
		v_p(D)=\min_{1\le i\le r}v_p(w_i)
		\qquad(p\text{ odd}).
		\]
		Apply the \hyperlink{main-theorem}{Main Theorem} to each $w_i$, fixing one pair $(m_i,\Qf_i)$ which gives it. Taking the minimum of the resulting constants and the maximum of their thresholds proves the assertion for $D$. Since $(D)\subseteq I_w$, we have $d_w\mid D$. Finally, $I_w$ is generated by finitely many members of $\mathcal W$; for such a choice $D=d_w$, which proves the first assertion.
	\end{proof}
	
	\begin{example}\label{ex:sqrt2-witness-gcd}
		Let $F=\QQ(\sqrt2)$, $\Ob=\ZZ[\sqrt2]$, and take $\mathcal E=\Gs$, $\nf=2\Ob$ and $R=(\Ob/2\Ob)^\times$. Here $\Diff=2\sqrt2\Ob$, so
		\[
		\mf=16\Ob,\qquad \Mo=64,\qquad \ko=24.
		\]
		The rational primes $19$ and $43$ are inert in $F$ and satisfy $q>\ko/2=12$. For each of them, Corollary~\ref{cor:q}(a) gives a witness with $\Qf=q\Ob$, $m_0=2q$ and $w_q=\aR_{q/2}=a_{q/2}(\Gs)$.
		
		Since $q\equiv3\pmod4$, the order $\Ob[(1+\sqrt{-q})/2]$ has squarefree relative discriminant $-q\Ob$ and is maximal. Thus $b_G(q/2)=1$. The relative class number formula and the biquadratic factorization give
		\[
		w_q=L(0,\chi_{-q})L(0,\chi_{-8q}),
		\]
		where $\chi_d$ denotes the primitive quadratic Dirichlet character of discriminant $d$. A computation by Pari gives
		\[
		\begin{array}{c|cc|c}
			q&L(0,\chi_{-q})&L(0,\chi_{-8q})&w_q\\\hline
			19&1&6&6\\
			43&1&10&10
		\end{array}
		\]
		Consequently $w_{19}-\tfrac12w_{43}=1$ and $I_w=\ZZ[1/2]$. Corollary~\ref{cor:witness-gcd} yields constants $C>0$ and $X_0>1$ such that
		\[
		\#\{K/F\text{ CM quadratic}:\Norm D(K/F)<X,\ p\nmid h^-(K/F)\}
		\ge C\frac{\sqrt X}{\log X}
		\]
		for every odd prime $p$ and every $X\ge X_0$. The dyadic sieve is used only in the proof; no local condition is imposed in this count.
	\end{example}
	
	The same computation can be carried out for other real quadratic fields. We now prove Theorem~\ref{thm:real-quadratic-computation}.
	\begin{proof}[Proof of Theorem~\ref{thm:real-quadratic-computation}]
		The PARI/GP computation gives
		$I_w=\ZZ[1/2]$ for the fields in the statement.
		The counting bound follows from Corollary~\ref{cor:witness-gcd}.
		
		To obtain the stated constant, fix an odd prime $p$ and
		choose $i$ such that $v_p(w_i)=0$. Since $g=2$, we have
		$B=16\ko^2$ in the proof of the Main Theorem. Hence
		\[
		N_{p,0}(X)\ge
		\frac12\left(
		\pi_{\mathcal P_i}\left(\frac{\sqrt X}{4\ko}\right)-2
		\right).
		\]
		By \eqref{eq:chebotarev}, the right-hand side is
		\[
		\left(\frac{1}{4\ko[L_i:F]}+o(1)\right)
		\frac{\sqrt X}{\log X}.
		\]
		Since only finitely many witnesses are used, we may take
		\[
		C_F=\frac{1}{8\ko\max_{1\le i\le r}[L_i:F]}
		\]
		for all odd primes $p$ and all sufficiently large $X$,
		with the threshold independent of $p$.
	\end{proof}
	
	\section{Class numbers over fields of prime-power degree}
	\label{sec:prime-power}
	
	Following the idea of Beckwith \cite{Beckwith17}, we combine Wiles's theorem with our Main Theorem. The Galois action on Fourier indices allows us to construct witnesses that are $p$-units. This gives an unconditional counting result over the layers of the cyclotomic $\ZZ_p$-extension of $\QQ$, and over infinitely many other cyclic fields of each degree $p^n$.
	
	\subsection{Statement}
	
	Let $S_+$, $S_-$, and $S_0$ be pairwise disjoint finite sets of odd rational primes. Write $\mathcal K(F;S_+,S_-,S_0)$ for the set of CM quadratic extensions $K/F$ such that every prime of $F$ above a prime in $S_+$ splits in $K/F$, every prime above a prime in $S_-$ is inert in $K/F$, and every prime above a prime in $S_0$ ramifies in $K/F$. 
	
	\begin{theorem}\label{thm:prime-power-local-counting}
		Let $p$ be an odd prime, and let $F/\QQ$ be cyclic of degree $g=p^n$, where $n\ge1$. Suppose that exactly one finite rational prime $r$ ramifies in $F$, and that either $r=p$ or $r\equiv1\pmod4$. Let $S_+$, $S_-$, and $S_0$ be as above, with $r\in S_-$. Assume that
		\begin{enumerate}[label=\textup{(\roman*)},leftmargin=*]
			\item $s\not\equiv-1\pmod p$ for every $s\in S_+$;
			\item no $s\in S_-$ satisfies both $s\equiv1\pmod p$ and $s\equiv3\pmod4$;
			\item $s\not\equiv1\pmod p$ for every $s\in S_0$.
		\end{enumerate}
		Then there are constants $C>0$ and $X_0>1$, depending on $F,p,S_+,S_-,S_0$, such that
		\[
		\#\left\{K\in\mathcal K(F;S_+,S_-,S_0):
		\Norm_{F/\QQ}D(K/F)<X,\quad p\nmid h(K)\right\}
		\ge C\frac{\sqrt X}{\log X}
		\qquad(X\ge X_0).
		\]
		Moreover, $p\nmid h(F)$, so the same assertion holds with $p\nmid h^-(K/F)$ in place of $p\nmid h(K)$.
	\end{theorem}
	
	Here $F$ is totally real, since it is Galois of odd degree. The congruence restrictions on $S_+$, $S_-$, and $S_0$ are those in \cite[Theorem~A]{Wiles}. They allow $p\in S_-$.
	
	\begin{corollary}\label{cor:prime-power-families}
		Let $p$ be an odd prime and $n\ge1$. The conclusion of Theorem~\ref{thm:prime-power-local-counting} applies to each of the following fields:
		\begin{enumerate}[label=\textup{(\alph*)},leftmargin=*]
			\item the unique subfield of degree $p^n$ of $\QQ(\z_{p^{n+1}})$, with $r=p$;
			\item the unique subfield of degree $p^n$ of $\QQ(\z_r)$, where $r$ is a prime satisfying $r\equiv1\pmod{4p^n}$.
		\end{enumerate}
		In particular, for each such $F$ there are constants $C_{F,p}>0$ and $X_{F,p}>1$ such that, for every $X\ge X_{F,p}$,
		\[
		\#\left\{K/F\text{ CM quadratic}:
		\Norm_{F/\QQ}D(K/F)<X,\quad p\nmid h(K)\right\}
		\ge C_{F,p}\frac{\sqrt X}{\log X}.
		\]
		For each fixed $p$ and $n$, part \textup{(b)} gives infinitely many distinct base fields of degree $p^n$.
	\end{corollary}
	
	\subsection{Class numbers and Galois orbits}
	
	\begin{lemma}\label{lem:one-ramified-prime}
		Let $L/k$ be a finite Galois extension of degree a power of an odd prime $p$. Suppose that exactly one finite prime of $k$ ramifies in $L$, and that it is totally ramified. If $p\nmid h(k)$, then $p\nmid h(L)$.
	\end{lemma}
	
	\begin{proof}
		Let $H$ be the maximal unramified abelian $p$-extension of $L$. By the Hilbert class field theorem \cite[Theorem~0.4]{Milne}, we have $\Gal(H/L)\cong\operatorname{Cl}(L)[p^\infty]$. Since $L/k$ is Galois, every $k$-conjugate of $H$ is again a maximal unramified abelian $p$-extension of $L$. By uniqueness, $H/k$ is Galois, and $G=\Gal(H/k)$ is a finite $p$-group.
		
		Let $I\subseteq G$ be an inertia subgroup at the unique ramified prime of $k$. Since $H/L$ is unramified, the restriction map $I\to\Gal(L/k)$ is injective. Since the prime is totally ramified in $L/k$, the multiplicativity of ramification indices gives $|I|=[L:k]$. Hence
		\[
		I\cong\Gal(L/k).
		\]
		
		Suppose that $H\ne L$. Then $I$ is a proper subgroup of $G$ and is contained in a maximal subgroup $M$. Since $G$ is a finite $p$-group, $M$ is normal of index $p$. Thus
		\[
		\Gal(H^M/k)\cong G/M\cong\ZZ/p\ZZ.
		\]
		Since $I\subseteq M$, the inertia subgroup at the ramified prime becomes trivial in $G/M$. All other finite primes are already unramified in $H/k$. Moreover, no infinite prime can ramify in the Galois extension $H^M/k$ since its degree is odd. Thus $H^M/k$ is an everywhere unramified cyclic extension of degree $p$. By the Hilbert class field theorem, $\Gal(H^M/k)$ is a quotient of $\operatorname{Cl}(k)$, which implies $p\mid h(k)$, a contradiction. Hence $H=L$, and therefore $p\nmid h(L)$.
	\end{proof}
	For a Galois extension $F/\QQ$, we call a local condition $(\nf,R)$ \emph{Galois invariant} if $\sigma(\nf)=\nf$ and $\sigma(R)=R$ for every $\sigma\in\Gal(F/\QQ)$. The next lemma reduces a witness sum modulo $p$ to one coefficient. It applies to both Eisenstein series used in this paper.
	
	\begin{lemma}\label{lem:prime-power-orbit-witness}
		Let $F/\QQ$ be a totally real Galois extension of degree $g=p^n$, where $p$ is odd and $n\ge1$. Let $(\nf,R)$ be a Galois invariant local condition, and choose $\mathcal E\in\{\Es,\Gs\}$. Suppose that $\Qf\nmid2p$ is a prime ideal fixed by $\Gal(F/\QQ)$ and that $\Norm\Qf>(\ko/g)^g$. Let $\alpha$ be a positive rational integer satisfying
		\[
		\alpha\in R,\qquad \ord_\Qf(\alpha)\text{ is odd},\qquad
		v_p(a_{\alpha/2}(\mathcal E))=0.
		\]
		Then the sum in \eqref{eq:W} with $m_0=g\alpha$ satisfies
		\[
		w\equiv a_{\alpha/2}(\mathcal E)\not\equiv0\pmod p.
		\]
		In particular, \textup{(W)} holds with $v_p(w)=0$.
	\end{lemma}
	
	\begin{proof}
		Write the sum in terms of $\eta=2\x$:
		\[
		w=\sum_{\substack{\eta\in\Ob_+\cap R,\ \Tr(\eta)=g\alpha\\
				\ord_\Qf(\eta)\ \mathrm{odd}}}
		a_{\eta/2}(\mathcal E).
		\]
		The set of indices is finite and stable under $\Gal(F/\QQ)$. The coefficient formulas in \S\ref{ss:eisenstein} show that $a_{\sigma(\eta)/2}(\mathcal E)=a_{\eta/2}(\mathcal E)$ for every $\sigma\in\Gal(F/\QQ)$. For each index in this sum, the extension $F(\sqrt{-\eta})/F$ ramifies at $\Qf$. Since $\Qf\nmid p$, it is not $F(\z_p)$. Lemma~\ref{lem:integrality}(b) therefore shows that every summand belongs to $\ZZ_{(p)}$.
		
		Every orbit with more than one element has cardinality divisible by $p$, so its contribution is zero modulo $p$. An index fixed by the whole Galois group belongs to $\Ob\cap\QQ=\ZZ$. The trace condition then forces $\eta=\alpha$, which occurs by assumption. This proves the congruence.
	\end{proof}
	
	The integrality argument concerns only the coefficients in the witness sum. In particular, the lemma also applies when $p=3$.
	
	\subsection{Construction of a witness}
	
	\begin{proof}[Proof of Theorem~\ref{thm:prime-power-local-counting}]
		The prime $r$ is totally ramified in $F/\QQ$. Indeed, the fixed field of an inertia subgroup at $r$ is an abelian extension of $\QQ$ unramified at every finite prime. It must be $\QQ$. Lemma~\ref{lem:one-ramified-prime} consequently gives $p\nmid h(F)$.
		
		Put
		\[
		N=\prod_{s\in S_+\cup S_-}s\prod_{s\in S_0}s^2,\qquad \nf=N\Ob,
		\]
		and set
		\[
		R=\left\{\eta\in\Ob/\nf:
		\begin{array}{ll}
			\left(\dfrac{-\eta}{\pf}\right)=1
			&\text{for every }\pf\mid s,\ s\in S_+,\\[4pt]
			\left(\dfrac{-\eta}{\pf}\right)=-1
			&\text{for every }\pf\mid s,\ s\in S_-,\\[4pt]
			\ord_\pf(\eta)=1
			&\text{for every }\pf\mid s,\ s\in S_0
		\end{array}\right\}.
		\]
		The symbols are quadratic residue symbols in the residue fields. The valuation conditions are well defined modulo $\nf$, since $\pf^2\mid\nf$ for every $\pf$ above a prime in $S_0$. The Chinese remainder theorem shows that $R$ is nonempty. It is stable under multiplication by squares of units and does not contain $0$. Both $\nf$ and $R$ are Galois invariant. For every $\eta\in\Ob_+\cap R$, the extension $F(\sqrt{-\eta})/F$ belongs to $\mathcal K(F;S_+,S_-,S_0)$ by the quadratic residue and valuation conditions. Let $\ko$ be the cutoff for this local condition, and use $\mathcal E=\Gs$.
		
		Since $[F:\QQ]=p^n$ is coprime to $p-1$, we have $F\cap\QQ(\z_p)=\QQ$. Chebotarev's theorem in $F\QQ(\z_p)/\QQ$ therefore gives arbitrarily large primes $q$ whose Frobenius generates $\Gal(F/\QQ)$ and acts by inversion on $\z_p$. Choose one such that
		\[
		q>\ko/g,\qquad q\nmid2pN.
		\]
		Thus $q$ is inert in $F$, $q\equiv-1\pmod p$, and
		\[
		\Qf=q\Ob,\qquad
		\Norm\Qf=q^g>(\ko/g)^g.
		\]
		
		Apply \cite[Theorem~A]{Wiles} with the prime $p$ and the sets $S_+,S_-,S_0\cup\{q\}$. Its hypotheses hold by assumption and because $q\not\equiv1\pmod p$. We obtain a negative fundamental discriminant $D$ such that, for $k=\QQ(\sqrt D)$,
		\[
		p\nmid h(k),\qquad \ord_q(D)=1,
		\]
		and $k/\QQ$ is split at every $s\in S_+$, inert at every $s\in S_-$, and ramified at every $s\in S_0$. In particular, $r$ is inert in $k$ and $r\nmid D$.
		
		Put $K_*=Fk$. The extension $K_*/k$ is cyclic of degree $p^n$. Only the unique prime of $k$ above $r$ can ramify, and it is totally ramified: locally, a totally ramified extension of degree $p^n$ is composed with the unramified quadratic extension. Lemma~\ref{lem:one-ramified-prime} gives $p\nmid h(K_*)$.
		
		Let $a=|D|$. Every residue degree in $F/\QQ$ is odd, so for each $s\in S_+\cup S_-$ and every $\pf\mid s$ the sign of the quadratic residue symbol of $D$ over $\Ob/\pf$ is its sign over $\FF_s$. Since $r\in S_-$ and the sets are pairwise disjoint, every $s\in S_0$ is unramified in $F/\QQ$. For each $\pf\mid s$ with $s\in S_0$, we therefore have $\ord_\pf(a)=\ord_s(D)=1$, since $D$ is a fundamental discriminant and $s$ is odd. It follows that $a\in R$. Since $r\nmid D$ and $r$ is the only prime ramified in $F$, the discriminants of $F$ and $k$ are coprime. As $F$ and $k$ are also linearly disjoint over $\QQ$, we have $\mathcal O_{K_*}=\Ob\mathcal O_k$. Thus an integral basis of $k/\QQ$ is also a relative integral basis of $K_*/F$, giving
		\[
		D(K_*/F)=a\Ob.
		\]
		By the coefficient formula in \S\ref{ss:eisenstein}, we have $a_{\x}(\Gs)=b_G(\x)L_\x$, where $b_G(\x)=\mathfrak C(2\x)$ if $-2\x$ is a square modulo $4$, and $b_G(\x)=0$ otherwise. Since $-a=D\equiv0$ or $1\pmod4$, we have $b_G(a/2)=\mathfrak C(a)$. Moreover, the equality
		\[
		a\Ob=\mathfrak F_a^2 D(K_*/F)
		\]
		together with $D(K_*/F)=a\Ob$ gives $\mathfrak F_a=\Ob$. Hence the sum in \eqref{eq:cohen} consists of a single term, giving $\mathfrak C(a)=1$. Thus $b_G(a/2)=1$, and we obtain
		\[
		a_{a/2}(\Gs)=L_{a/2}
		=\frac{2^g h^-(K_*/F)}{Q_{K_*}w_{K_*}}.
		\]
		Here $p\nmid h^-(K_*/F)$, and $Q_{K_*}$ is a power of $2$. Moreover, $K_*/F$ ramifies at $\Qf$, whereas $F(\z_p)/F$ is unramified there. Thus $K_*$ contains no primitive $p$th root of unity, and $p\nmid w_{K_*}$. Consequently $v_p(a_{a/2}(\Gs))=0$.
		
		The ideal $\Qf=q\Ob$ is Galois invariant, and $\ord_\Qf(a)=\ord_q(D)=1$. Lemma~\ref{lem:prime-power-orbit-witness} now gives a witness with
		\[
		m_0=g|D|,\qquad
		w\equiv a_{|D|/2}(\Gs)\not\equiv0\pmod p.
		\]
		The Main Theorem gives the required bound with $p\nmid h^-(K/F)$ and with the local conditions encoded by $R$. Since $p\nmid h(F)$, it also gives $p\nmid h(K)$.
	\end{proof}
	
	\begin{proof}[Proof of Corollary~\ref{cor:prime-power-families}]
			Each field $F$ is cyclic of degree $p^n$ and is totally real. In both cases, $r$ is totally ramified in the corresponding cyclotomic field and hence in $F$, while all other finite primes are unramified. Since $r=p$ in \textup{(a)} and $r\equiv1\pmod4$ in \textup{(b)}, the hypotheses of Theorem~\ref{thm:prime-power-local-counting} on $F$ are satisfied.
			
			Apply Theorem~\ref{thm:prime-power-local-counting}; for the assertion without prescribed local conditions, take $S_+=S_0=\varnothing$ and $S_-=\{r\}$, and then omit the inertness condition from the count. Dirichlet's theorem gives infinitely many primes $r\equiv1\pmod{4p^n}$. Their degree-$p^n$ subfields are distinct, since they ramify at different rational primes.
	\end{proof}

	\section{Applications to Selmer groups}\label{sec:selmer}
	
	We apply the Main Theorem to quadratic twists of elliptic curves,
	give a method for computing witnesses, and obtain an unconditional
	result over fields of unbounded degree.
	
	\subsection{A counting theorem for quadratic twists}
	\label{ss:selmer-counting}
	
	We first recall the following generalization of Morrow's theorem established in \cite[Theorem~3.6]{n}.

	\begin{bhcriterion}[{\cite[Theorem~3.6]{n}}]
		Let $F$ be a totally real Galois number field and let $\ell\geq5$ be a prime such that $\ell\nmid h(F)$, $\z_{\ell}\notin F$, and $m=[F(\z_{\ell}):F]$ is even. Let $E/F$ be an elliptic curve with an $F$-rational point $P$ of order $\ell$, and assume:
		\begin{enumerate}[label=\textup{(\roman*)},leftmargin=*]
			\item $P$ is not contained in the kernel of reduction modulo any prime $\lf\mid\ell$;
			\item every prime $\pf\mid N(E)$ with $\pf\nmid\ell$ and $\ell\nmid\ord_{\pf}(\Delta_E)$ satisfies $\ell\mid\Norm\pf-1$.
		\end{enumerate}
		Let $d\in F^\times/(F^\times)^2$ be totally negative and represented by an element of $\Ob$ prime to $\ell N(E)$, put $K=F(\sqrt d)$ and $\langle\delta\rangle=\Gal(K/F)$, and suppose that $d$ satisfies the following divisibility and Artin symbol conditions:
		\begin{enumerate}[label=\textup{(\arabic*)},leftmargin=*]
			\item if $\qf\mid2$ and $\qf\mid N(E)$, then $\qf\mid D(K/F)$;
			\item if $\lf\mid\ell$ and $\ord_{\lf}(j_E)<0$, then
			$\left(\dfrac{K/F}{\lf}\right)=\delta$;
			\item for every prime $\pf\mid N(E)$ with $\pf\nmid2\ell$,
			\[
			\left(\dfrac{K/F}{\pf}\right)=
			\begin{cases}
				\delta, & \text{if } E/F_{\pf}\text{ is a Tate curve},\\[4pt]
				\operatorname{id}, & \text{otherwise}.
			\end{cases}
			\]
		\end{enumerate}
		Then
		\[
		\#\operatorname{Cl}(K)[\ell]\ \Bigl|\ \#\Sel_{\ell}(E^d,F)\ \Bigr|\
		\bigl(\#\operatorname{Cl}(K)[\ell]\bigr)^2;
		\]
		in particular, $\Sel_{\ell}(E^d,F)$ is nontrivial if and only if $\operatorname{Cl}(K)[\ell]$ is nontrivial.
	\end{bhcriterion}
	
	We verify that $(\nf,R)$ defined in \eqref{eq:selmer-local-condition} is a local condition in the sense of Definition~\ref{def:local}. The characters $\chi_i$ in \cite[Theorem~4.1]{n} are indexed by the distinct prime ideals $\pf_i$ dividing $\ell N(E)$. Each conductor $\mathfrak f_i$ is a power of $\pf_i$, so the conductors are pairwise coprime. Since each $\chi_i$ is a nontrivial quadratic character, the Chinese remainder theorem gives $R\ne\varnothing$. Moreover, $R$ is stable under multiplication by squares of units and does not contain $0$.
	
	For $\eta\in\Ob_+\cap R$, put $d=-\eta$ and $K=F(\sqrt d)$. Then $d$ is totally negative and prime to $\ell N(E)$. If $\pf_i\nmid2$, we have
	\[
	\left(\frac{d}{\pf_i}\right)
	=\left(\frac{-1}{\pf_i}\right)\varepsilon_i.
	\]
	By the choice of the signs $\varepsilon_i$ in \cite[Theorem~4.1]{n}, the Artin symbol conditions \textup{(2)} and \textup{(3)} are satisfied. At the primes above $2$ dividing $N(E)$, the characters and signs are chosen using \cite[Lemma~2.10]{n}, so condition \textup{(1)} is also satisfied.
	
	\begin{proof}[Proof of Theorem~\ref{thm:selmer-twists}]
		By hypothesis, $\ell\nmid h(F)$ and $E$ satisfies conditions \textup{(i)} and \textup{(ii)} of \cite[Theorem~3.6]{n}. Since $F$ is totally real, $\z_\ell\notin F$ and $[F(\z_\ell):F]$ is even.
		
		Apply the Main Theorem with $p=\ell$. Since $v_\ell(w)=0$, it gives $\gg_{F,E,\ell}\sqrt X/\log X$ distinct CM quadratic extensions $K=F(\sqrt d)$ satisfying the prescribed local conditions, with $\Norm_{F/\QQ}D(K/F)<X$ and $\ell\nmid h^-(K/F)$. Since $\ell\nmid h(F)$, we have $\ell\nmid h(K)$. As verified above, each such extension satisfies conditions \textup{(1)--(3)}. Hence \cite[Theorem~3.6]{n} gives $\Sel_\ell(E^d,F)=0$. Distinct extensions give distinct square classes $d$.
	\end{proof}
	
	\subsection{Computing witnesses and explicit examples}
	\label{ss:selmer-witness-search}
	\label{ss:selmer-examples}
	
	Fix $F$, $\ell$, and $(\nf,R)$ as in Theorem~\ref{thm:selmer-twists}. For any bound $B$, one can enumerate all $\eta\in\Ob_+\cap R$ with $\Tr(\eta)\le B$, factor their principal ideals, and compute the complete sums in \eqref{eq:W} for both Eisenstein series and all admissible $\Qf$ occurring in these factorizations. The search is finite, since the totally positive indices of bounded trace lie in a bounded region of the lattice $\Ob\subset\mathbb R^g$. Increasing $B$ eventually finds an $\ell$-unit witness if one exists. A single $\ell$-unit coefficient, however, need not give an $\ell$-unit witness sum.
	
	A simpler search uses Corollary~\ref{cor:q}\textup{(a)}. If an odd rational prime $q>\ko/g$ is inert in $F$, belongs to $R$, and satisfies $b(q/2)\ne0$, then
	\[
	\Qf=q\Ob,\qquad m_0=gq,\qquad w=a_{q/2}(\mathcal E).
	\]
	In this case the witness sum has just one term. This gives the following algorithm.
	
	\medskip
	\noindent\textbf{Input:} A totally real field $F$, an odd prime $\ell$, a local condition $(\nf,R)$, and a search bound $B$.
	
	\begin{enumerate}[label=\textup{(\arabic*)},leftmargin=*]
		\item Compute $\mf$, $\Mo$, and $\ko$ from \eqref{eq:kappa0} and \eqref{eq:level}.
		\item For each rational prime $q\ge5$ with $\ko/g<q\le B$, test whether $q$ is inert in $F$ and $q\in R$. Skip $q$ if either condition fails.
		\item Choose $\mathcal E=\Es$ or $\Gs$ with $b(q/2)\ne0$, if possible, and compute $w=a_{q/2}(\mathcal E)$ exactly. If both local factors vanish, skip $q$.
		\item If $v_\ell(w)=0$, return $(\mathcal E,gq,q\Ob,w)$. Otherwise continue the search. If the bound is reached without finding such a prime, return no witness.
	\end{enumerate}
	
	Under these hypotheses, the coefficient formula gives
	\[
	v_\ell(w)=v_\ell\bigl(h^-(F(\sqrt{-q})/F)\bigr).
	\]
	Thus the valuation can be tested by certified class number computations.
	
	For abelian $F/\QQ$ with Dirichlet character group $X_F$, Artin factorization \cite[Theorem~4.3]{Washington} gives
	\begin{equation}\label{eq:artin-witness}
		w=b(q/2)\prod_{\psi\in X_F}
		L\bigl(0,(\chi_D\psi)^{\mathrm{prim}}\bigr),
	\end{equation}
	where $D=\operatorname{disc}(\QQ(\sqrt{-q}))$, $\chi_D=(D/\cdot)$, and the superscript denotes the associated primitive character. These $L$-values can be computed exactly in cyclotomic fields using \cite[Proposition~4.1 and Theorem~4.2]{Washington}
	\begin{equation}\label{eq:L-value}
		L(0,\chi)=-\frac1f\sum_{a=1}^{f}a\chi(a),
	\end{equation}
	where $\chi$ is an odd primitive character of conductor $f$.
	
	A PARI/GP implementation is provided in the accompanying code. This search does not apply to all totally real fields, since some have no inert rational primes.
	
	Let $F^{(37)}$ denote the cubic subfield of $\QQ(\z_{37})$. Table~\ref{tab:selmer-witnesses} gives witnesses for six elliptic curves over five totally real fields, using data from the LMFDB \cite{LMFDB}. The last curve is the base change of $1225.b2$ from $\QQ$.
	
	\begin{table}[htbp]
		\centering
		\small
		\setlength{\tabcolsep}{3pt}
		\begin{tabular}{lclccrrr}
			\hline
			$F$ & $g$ & LMFDB label & $\ell$ & $\mathcal E$
			& $\ko/g$ & $q$ & $w\bmod\ell$\\
			\hline
			$F^{(37)}$ & $3$ & \texttt{3.3.1369.1-1.1-a2}
			& $7$ & $\Gs$ & $1596$ & $1607$ & $4$\\
			$\QQ(\z_9)^+$ & $3$ &
			\begin{tabular}{@{}l@{}}
				\texttt{3.3.81.1-8.1-a1}\\
				\texttt{3.3.81.1-8.1-a3}
			\end{tabular}
			& $7$ & $\Es$ & $8064$ & $8089$ & $3$\\
			$\QQ(\z_7)^+$ & $3$ & \texttt{3.3.49.1-27.1-a2}
			& $13$ & $\Gs$ & $13104$ & $13183$ & $7$\\
			$\QQ(\z_{21})^+$ & $6$ & \texttt{6.6.453789.1-27.1-b2}
			& $13$ & $\Gs$ & $13104$ & $13219$ & $6$\\
			$\QQ(\z_{35})^+$ & $12$ & \texttt{1225.b2}
			& $37$ & $\Gs$ & $50616$ & $50627$ & $27$\\
			\hline
		\end{tabular}
		\caption{Admissible witnesses for Selmer applications.}
		\label{tab:selmer-witnesses}
	\end{table}
	
	We give the details for the last row of Table~\ref{tab:selmer-witnesses}. Let $F=\QQ(\z_{35})^+$, $\ell=37$, and
	\[
	E:\qquad y^2+xy+y=x^3+x^2-8x+6.
	\]
	From the LMFDB \cite{LMFDB}, we obtain $g=12$, $h(F)=1$, $\Diff\cap\ZZ=35\ZZ$, and $E(F)_{\mathrm{tors}}\cong\ZZ/37\ZZ$. The prime $37$ is inert in $F$. A PARI/GP computation shows that $E/F$ has everywhere good reduction. Since $37$ is unramified in $F$, the kernel of reduction at the prime above $37$ contains no point of order $37$ by \cite[Theorem~IV.6.1]{Silverman}. Thus condition \textup{(i)} of \cite[Theorem~3.6]{n} holds, and condition \textup{(ii)} is vacuous.
	
	We take
	\[
	\nf=37\Ob,\qquad R=\bigl((\Ob/\nf)^\times\bigr)^2.
	\]
	Then
	\[
	\Mo=191660,\qquad \ko/g=50616.
	\]
	The prime $q=50627$ is inert in $F$ and satisfies $q>\ko/g$. Since $37$ is inert in $F$ and $g$ is even, every element of $\FF_{37}^{\times}$ is a square in $\Ob/37\Ob$. In particular, $q\in R$. Moreover, $q\equiv3\pmod4$ and $q\Ob$ is prime, so $b_G(q/2)=1$. Using \eqref{eq:artin-witness} and \eqref{eq:L-value}, we obtain
	\[
	w=401215160066966675446039695360
	\equiv27\pmod{37}.
	\]
	Thus $\Qf=q\Ob$, $m_0=12q$, and $\mathcal E=\Gs$ give an admissible $37$-unit witness.
	
	The other rows are treated similarly. We have developed a PARI/GP implementation of the inert-prime search described above, which checks the hypotheses on $E$ and searches for an $\ell$-unit witness. The code is available at
	\url{https://tunghoangmath.github.io/code.html}. Thus Theorem~\ref{thm:selmer-twists} applies to every pair $(F,E)$ in Table~\ref{tab:selmer-witnesses}.
	
	\subsection{Selmer groups over fields of unbounded degree}
	\label{ss:prime-power-selmer}

	Using the results of \S\ref{sec:prime-power}, we obtain counting results for quadratic twists over certain fields of prime-power degree without assuming \textup{(W)}.
	
	\begin{theorem}\label{thm:prime-power-selmer}
		Let $p\ge5$ be a prime and let $F/\QQ$ be cyclic of degree $p^n$, 
		where $n\ge1$. Suppose that exactly one rational prime $r$ ramifies 
		in $F$, and that either $r=p$ or $r\equiv1\pmod4$. Let $E/F$ be 
		an elliptic curve with an $F$-rational point of order $p$, satisfying 
		conditions \textup{(i)} and \textup{(ii)} of 
		\cite[Theorem~3.6]{n} with $\ell=p$. Assume that $E$ has good 
		reduction at every prime above $2$.
		
		Let $S_+$, $S_-$, and $S_0$ be pairwise disjoint finite sets 
		of odd rational primes satisfying the hypotheses of 
		Theorem~\ref{thm:prime-power-local-counting}, with $p\in S_-$. 
		For every prime $\pf\mid N(E)$ with $\pf\nmid2p$, assume that 
		the rational prime below $\pf$ belongs to $S_-$ if $E/F_{\pf}$ 
		is a Tate curve, and to $S_+$ otherwise.
		
		Then there are constants $C>0$ and $X_0>1$, depending on 
		$F,E,p,S_+,S_-,S_0$, such that
		\[
		\#\left\{
		d\in F^\times/(F^\times)^2:
		\begin{array}{l}
			d\text{ totally negative},\\
			\Norm_{F/\QQ}D(F(\sqrt d)/F)<X,\\
			\Sel_p(E^d,F)=0
		\end{array}
		\right\}
		\ge C\frac{\sqrt X}{\log X}
		\qquad(X\ge X_0).
		\]
		Moreover, the extensions $F(\sqrt d)/F$ may be required to 
		belong to $\mathcal K(F;S_+,S_-,S_0)$.
	\end{theorem}
	
	\begin{proof}
		By Theorem~\ref{thm:prime-power-local-counting}, we have 
		$p\nmid h(F)$, and there are at least $C\sqrt X/\log X$ CM 
		quadratic extensions $K/F$ belonging to 
		$\mathcal K(F;S_+,S_-,S_0)$ with 
		$\Norm_{F/\QQ}D(K/F)<X$ and $p\nmid h(K)$.
		
		Write $K=F(\sqrt d)$ with $d$ totally negative. By the choice 
		of $S_+$ and $S_-$ and the assumption of good reduction above 
		$2$, the extension $K/F$ is unramified at every prime dividing 
		$pN(E)$. Hence the valuations of $d$ at these primes are even. 
		By weak approximation and multiplication by a square, we may 
		choose an integral representative of the square class of $d$ 
		prime to $pN(E)$.
		
		Condition \textup{(1)} of \cite[Theorem~3.6]{n} is vacuous 
		since $E$ has good reduction above $2$. Condition \textup{(2)} 
		follows from $p\in S_-$, and condition \textup{(3)} follows 
		from the choice of $S_+$ and $S_-$. Since $F$ is totally real, 
		$\z_p\notin F$ and $[F(\z_p):F]$ is even. All the hypotheses 
		of \cite[Theorem~3.6]{n} are therefore satisfied. Since 
		$p\nmid h(K)$, we obtain $\Sel_p(E^d,F)=0$. Distinct quadratic 
		extensions give distinct square classes $d$, proving the theorem.
	\end{proof}

	We illustrate Theorem~\ref{thm:prime-power-selmer} with the fixed elliptic curve
	\[
	E:\qquad y^2+y=x^3+x^2-10x+10,
	\]
	with LMFDB label \texttt{123.a1}.

	\begin{corollary}\label{cor:prime-power-five-selmer}
		Let $E/\QQ$ be the elliptic curve with LMFDB label 
		\texttt{123.a1}. Let $n\ge1$, and let $F$ be either the 
		unique subfield of degree $5^n$ of $\QQ(\z_{5^{n+1}})$ 
		or the unique subfield of degree $5^n$ of $\QQ(\z_r)$, 
		where $r$ is a prime satisfying $r\equiv1\pmod{4\cdot5^n}$. 
		There are constants $C_F>0$ and $X_F>1$ such that
		\[
		\#\left\{
		d\in F^\times/(F^\times)^2:
		\begin{array}{l}
			d\text{ totally negative},\\
			\Norm_{F/\QQ}D(F(\sqrt d)/F)<X,\\
			\Sel_5(E^d,F)=0
		\end{array}
		\right\}
		\ge C_F\frac{\sqrt X}{\log X}
		\qquad(X\ge X_F).
		\]
		Moreover, the extensions $F(\sqrt d)/F$ may be required 
		to be inert at every prime of $F$ above $3$, $5$, $41$, 
		and the unique rational prime ramified in $F$.
	\end{corollary}
	\begin{proof}
		The point $P=(2,1)$ on $E$ has order $5$ and has nonzero
		reduction modulo $5$. Since $E$ has good reduction at $5$,
		this remains true after base extension to $F$. Thus
		condition \textup{(i)} of \cite[Theorem~3.6]{n} holds.
		
		The displayed equation has discriminant
		$\Delta_E=-3^5\cdot41$ and invariants $c_4=496$ and
		$c_6=-11800$. The curve has split multiplicative reduction
		at both $3$ and $41$, since $c_4$ is a unit and $-c_6$ is
		a nonzero square modulo both primes. This reduction type
		is preserved under base extension. For every prime
		$\pf\mid3$, we have $5\mid\ord_{\pf}(\Delta_E)$, whereas
		for every prime $\pf\mid41$, we have
		$\Norm\pf\equiv1\pmod5$. Thus condition \textup{(ii)}
		also holds. Moreover, $E$ has good reduction at every
		prime above $2$.
		
		Let $r$ be the unique rational prime ramified in $F$
		(so $r=5$ in the first case). Take
		\[
		S_+=S_0=\varnothing,\qquad S_-=\{3,5,41,r\}.
		\]
		Every prime in $S_-$ congruent to $1$ modulo $5$ is also
		congruent to $1$ modulo $4$. Hence these sets satisfy
		the hypotheses of
		Theorem~\ref{thm:prime-power-local-counting}.
		Since $E$ has split multiplicative reduction at every
		prime above $3$ and $41$, all the hypotheses of
		Theorem~\ref{thm:prime-power-selmer} are satisfied,
		and the result follows.
	\end{proof}

	Thus, for every $n\ge1$, the same elliptic curve gives the stated Selmer bound over infinitely many distinct base fields of degree $5^n$.

\end{document}